\documentclass[12pt]{article}
\usepackage{latexsym, amssymb, amsmath, amscd, amsfonts, epsfig, graphicx, colordvi,verbatim,ifpdf,extarrows}
\usepackage{amsfonts, amsmath, amssymb}
\usepackage{amssymb,amsfonts,amsmath,latexsym,epsfig,cite, psfrag,eepic,color}
\usepackage{amscd,graphics}
\usepackage{latexsym, amssymb,  amsmath,amscd, amsfonts, epsfig, graphicx, colordvi,amsthm}

\usepackage{graphicx}
\usepackage{epstopdf}
\usepackage{color}
\usepackage{ifpdf}
\usepackage{fancybox}
\usepackage[font=small,labelfont=bf,labelsep=none]{caption}
\usepackage{float}

\newtheorem{thm}{Theorem}[section]
\newtheorem{prop}[thm]{Proposition}
\newtheorem{conj}[thm]{Conjecture}

\newtheorem{defi}[thm]{Definition}
\newtheorem{lem}[thm]{Lemma}
\newtheorem{core}[thm]{Corollary}

\def\pf{\noindent{\it Proof.} }
\def\qed{\nopagebreak\hfill{\rule{4pt}{7pt}}
\medbreak}

\numberwithin{equation}{section}

\def\qed{\nopagebreak\hfill{\rule{4pt}{7pt}}
\medbreak}

\newlength{\boxedparwidth}
  {\begin{center} \begin{tabular}{|@{\hspace{.315in}}c@{\hspace{.15in}}|}
                  \hline \\ \begin{minipage}[t]{\boxedparwidth}
                  \setlength{\parindent}{.25in}}%
  {\end{minipage} \\ \\ \hline \end{tabular} \end{center}}

\begin{document}

\begin{center}

 {\Large \bf A bijection related to Bressoud's conjecture}
\end{center}

\begin{center}
 {Y.H. Chen}$^{1}$ and {Thomas Y. He}$^{2}$ \vskip 2mm

$^{1,2}$ School of Mathematical Sciences, Sichuan Normal University, Chengdu 610066, P.R. China

   \vskip 2mm

   $^1$chenyh@stu.sicnu.edu.cn, $^2$heyao@sicnu.edu.cn
\end{center}

{\noindent \bf Abstract.} Bressoud introduced the partition function $B(\alpha_1,\ldots,\alpha_\lambda;\eta,k,r;n)$, which counts the number of partitions with certain difference
conditions. Bressoud posed a conjecture on the generating function for the partition function $B(\alpha_1,\ldots,\alpha_\lambda;\eta,k,r;n)$ in multi-summation form.
In this article, we introduce a bijection related to Bressoud's conjecture. As an application, we give a new companion to the  G\"ollnitz-Gordon identities.



\section{Introduction}

A partition $\pi=(\pi_1,\pi_2,\ldots,\pi_\ell)$ of a positive integer $n$ is a finite non-increasing sequence of positive integers   such that $\pi_1+\pi_2+\cdots+\pi_\ell=n.$ The $\pi_i$ are called the parts of $\pi$. Let $\ell(\pi)$ be the number of parts of $\pi$ and let $|\pi|$ be the sum of parts of $\pi$.

Assume that $\alpha_1,\alpha_2,\ldots, \alpha_\lambda$ and $\eta$ are integers such that \[0<\alpha_1<\alpha_2<\cdots<\alpha_\lambda<\eta,\text{ and }\alpha_i=\eta-\alpha_{\lambda+1-i}\text{ for }1\leq i\leq \lambda.\]
Bressoud \cite{Bressoud-1980} introduced the partition function $B(\alpha_1,\ldots,\alpha_\lambda;\eta,k,r;n)$.

\begin{defi}[Bressoud]\label{Bressoud-b-defi-0000000001}  For $k\geq r\geq\lambda\geq0,$ define
$B(\alpha_1,\ldots,\alpha_\lambda;\eta,k,r;n)$  to be the number of partitions $\pi=(\pi_1,\pi_2,\ldots,\pi_\ell)$ of $n$   satisfying the following conditions{\rm{:}}
\begin{itemize}
\item[{\rm (1)}] For $1\leq i\leq \ell$, $\pi_i\equiv0,\alpha_1,\ldots,\alpha_\lambda\pmod{\eta}${\rm{;}}

\item[{\rm (2)}] Only multiples of $\eta$ may be repeated{\rm{;}}

\item[{\rm (3)}]  For $1\leq i\leq\ell-k+1$, $ \pi_i\geq\pi_{i+k-1}+\eta$  with strict inequality if $\eta\mid\pi_i${\rm{;}}

\item[{\rm (4)}] At most $r-1$ parts of the $\pi_i$ are less than or equal to $\eta${\rm{.}}

  \end{itemize}
 \end{defi}

 Bressoud \cite{Bressoud-1980}  posed the following conjecture.

\begin{conj}[Bressoud] \label{Bressoud-gen-b-e} For $k\geq r\geq\lambda\geq0,$
\begin{equation*}
 \begin{split}
 &\sum_{n\geq0}B(\alpha_1,\ldots,\alpha_\lambda;\eta,k,r;n)q^n\\
 &\qquad=\sum_{N_1\geq\cdots\geq N_{k-1}\geq0}\frac{q^{\eta(N_1^2+\cdots+N_{k-1}^2+N_r+\cdots+N_{k-1})}}{(q^\eta;q^\eta)_{N_1-N_2}\cdots(q^\eta;q^\eta)_{N_{k-2}-N_{k-1}}(q^{\eta};q^{\eta})_{N_{k-1}}}\\
 &\qquad\qquad\times\prod_{s=1}^{\lambda}(-q^{\eta-\alpha_s-\eta N_s};q^\eta)_{N_s}\prod_{s=2}^{\lambda}(-q^{\eta -\alpha_s+\eta N_{s-1}};q^\eta)_\infty.
 \end{split}
 \end{equation*}
\end{conj}

Here and in the sequel, we assume that $|q|<1$ and employ the standard notation\cite{Andrews-1976}:
\[(a;q)_\infty=\prod_{i=0}^{\infty}(1-aq^i), \quad (a;q)_n=\frac{(a;q)_\infty}{(aq^n;q)_\infty},\]
and
\[(a_1,a_2,\ldots,a_m;q)_\infty=(a_1;q)_\infty(a_2;q)_\infty\cdots(a_m;q)_\infty.
\]

Kim and Yee \cite{Kim-Yee-2014}  gave a proof of Conjecture \ref{Bressoud-gen-b-e}
  for   $\lambda=2$ with the aid of Gordon markings introduced by  Kur\c{s}ung\"oz   \cite{Kursungoz-2010a, Kursungoz-2010}.  Recently, Kim  \cite{Kim-2018} established  the following generating function for $B(\alpha_1,\ldots,\alpha_\lambda;\eta,k,r;n)$ in infinite product form.

  \begin{thm}[Kim]\label{bressoud-product}
  For $k\geq r\geq\lambda\geq0$,
\begin{equation*}
 \begin{split}
 &\sum_{n\geq0}B(\alpha_1,\ldots,\alpha_\lambda;\eta,k,r;n)q^n\\[5pt]
 &\quad=\frac{(-q^{\alpha_1},\ldots,-q^{\alpha_\lambda};q^{\eta})_\infty
 (q^{\eta(r-\frac{\lambda}{2})},q^{\eta(2k-r-\frac{\lambda}{2}+1)},
 q^{\eta(2k-\lambda+1)};q^{\eta(2k-\lambda+1)})_\infty}{(q^\eta;q^\eta)_\infty}.
 \end{split}
 \end{equation*}
 \end{thm}

Then, Conjecture \ref{Bressoud-gen-b-e}
 is an immediate consequence of Theorem \ref{bressoud-product} and the following theorem obtained by Bressoud \cite{Bressoud-1980}. It is worth mentioning that the following theorem can specialize to many well-known Rogers-Ramanujan type identities, such as Rogers-Ramanujan-Gordon identities \cite{Gordon-1961} and G\"ollnitz-Gordon identities \cite{Gollnitz-1960,Gollnitz-1967,Gordon-1962,Gordon-1965}.
  \begin{thm}[Bressoud]\label{bressoud-double-sum}
 For  $k\geq r\geq\lambda\geq0$,
\begin{equation*}\label{Bressoud-conj-e}
 \begin{split}
 &\sum_{N_1\geq\cdots \geq N_{k-1}\geq0}\frac{q^{\eta(N_1^2+\cdots+N_{k-1}^2+N_r
 +\cdots+N_{k-1})}}{(q^\eta;q^\eta)_{N_1-N_2}
 \cdots(q^\eta;q^\eta)_{N_{k-2}-N_{k-1}}(q^{\eta};q^{\eta})_{N_{k-1}}}\\[5pt]
 &\qquad\qquad\times\prod_{s=1}^{\lambda}(-q^{\eta-\alpha_s-\eta N_s};q^\eta)_{N_s}\prod_{s=2}^{\lambda}(-q^{\eta -\alpha_s+\eta N_{s-1}};q^\eta)_\infty\\[10pt]
  &\qquad\qquad\qquad=\frac{(-q^{\alpha_1},\ldots,-q^{\alpha_\lambda};q^{\eta})_\infty(q^{\eta(r-\frac{\lambda}{2})},q^{\eta
 (2k-r-\frac{\lambda}{2}+1)},q^{\eta(2k-\lambda+1)}
 ;q^{\eta(2k-\lambda+1)})_\infty}{(q^\eta;q^\eta)_\infty}.
 \end{split}
 \end{equation*}
 \end{thm}

However, Conjecture \ref{Bressoud-gen-b-e} also cries out for a direct combinatorial proof. To do this, for $N_1\geq N_2\geq\cdots\geq N_{k-1}\geq 0$, the main task is to merge the partitions whose generating functions are
\begin{equation}\label{intro-0-2}
(-q^{\eta-\alpha_s-\eta N_s};q^\eta)_{N_s},
\end{equation}
\begin{equation}\label{intro-0-3}
 (-q^{\eta -\alpha_s+\eta N_{s-1}};q^\eta)_\infty,
 \end{equation}
and
\begin{equation}\label{intro-1}
\frac{q^{\eta(N_1^2+\cdots+N_{k-1}^2+N_r+\cdots+N_{k-1})}}{(q^\eta;q^\eta)_{N_1-N_2}\cdots(q^\eta;q^\eta)_{N_{k-2}-N_{k-1}}(q^{\eta};q^{\eta})_{N_{k-1}}}.
\end{equation}

In \cite{Kim-Yee-2014}, Kim and Yee showed how to merge the partitions whose generating functions are \eqref{intro-0-2} with $s=1$, \eqref{intro-0-2} with $s=2$, \eqref{intro-0-3} with $s=2$ and \eqref{intro-1}. In this article, we will introduce a bijection which tell us how to merge the partitions with generating function $(-q^{\eta-\alpha_{3}+\eta N_{2}};q^\eta)_\infty$ (\eqref{intro-0-3} with $s=3$) and the partitions whose generating function is given in  \eqref{intro-1}. For easier expression, we investigate the following identities:

\begin{equation*}\label{main-ins-1}
(-q^{1+2N_2};q^2)_\infty,
\end{equation*}
and
\begin{equation*}\label{main-ins-2}
\frac{q^{2(N_1^2+\cdots+N_{k-1}^2+N_r+\cdots+N_{k-1})}}{(q^2;q^2)_{N_1-N_2}\cdots(q^2;q^2)_{N_{k-2}-N_{k-1}}(q^{2};q^{2})_{N_{k-1}}}.
\end{equation*}

The main result of this article is given below.
\begin{thm}\label{equiv-main-lem-1-1}
For $k\geq r \geq 3$ and $p,t\geq0,$ there is a bijection $\Phi_{p,t}$ between $\mathbb{C}_{<}(k,r|p,t)$ and $\mathbb{C}_{=}(k,r|p,t)$. Moreover, for a partition $\pi\in \mathbb{C}_{<}(k,r|p,t)$, we have $\omega=\Phi_{p,t}(\pi)\in\mathbb{C}_{=}(k,r|p,t)$ such that
\[|\omega|=|\pi|+2p+2t+1 \text{ and } \ell(\omega)=\ell(\pi)+1.\]
\end{thm}
Since the explicit definitions of $\mathbb{C}_{<}(k,r|p,t)$ and $\mathbb{C}_{=}(k,r|p,t)$ are complicate, we put them in Section 2. Let $B(1;2,3,3;\ell,n)$ denote the number of partitions counted by $B(1;2,3,3;n)$ with exactly $\ell$ parts. Based on the bijection $\Phi_{p,t}$ in Theorem \ref{equiv-main-lem-1-1}, we will give the following formula for the generating function of $B(1;2,3,3;\ell,n)$, which can be regarded as a new companion to the generalizations of the G\"ollnitz-Gordon identities \cite{Bressoud-1980}.

\begin{thm}\label{main-thm} The generating function of $B(1;2,3,3;\ell,n)$ is
\begin{equation*}\label{main-eqn}
\displaystyle\sum_{\ell,n\geq0}B(1;2,3,3;\ell,n)x^\ell q^n=\sum_{N_1\geq N_{2}\geq0}\frac{q^{2(N_1^2+N_{2}^2)}(-xq^{1+2N_2};q^2)_\infty x^{N_1+N_{2}}}{(q^2;q^2)_{N_1-N_2} (q^{2};q^{2})_{N_{2}}}.
\end{equation*}
\end{thm}

 This article is organized as follows. In Section 2, we give the explicit definitions of $\mathbb{C}_{<}(k,r|p,t)$ and $\mathbb{C}_{=}(k,r|p,t)$ in Theorem \ref{equiv-main-lem-1-1} and investigate the properties of them. Furthermore, we give equivalent statements of Theorem \ref{equiv-main-lem-1-1}, which are stated in Theorems \ref{equiv-main-H} and \ref{equiv-main-I}. Section 3 is devoted to introducing the dilation operation and the reduction operation, which allow us to provide a proof of Theorem \ref{equiv-main-H}. In Section 4, we introduce  the insertion operation and the separation operation and then give a proof of Theorem \ref{equiv-main-I}. Finally, we give a combinatorial proof of Theorem \ref{main-thm} in Section 5.

 \section{$\mathbb{C}_{<}(k,r|p,t)$ and $\mathbb{C}_{=}(k,r|p,t)$}
This section is devoted to giving the explicit definitions of  $\mathbb{C}_{<}(k,r|p,t)$ and $\mathbb{C}_{=}(k,r|p,t)$ in Theorem \ref{equiv-main-lem-1-1}. To do this, we need to recall the definition of G\"ollnitz-Gordon marking given in \cite{He-Zhao-2023} and introduce the starting types based on G\"ollnitz-Gordon marking.

\subsection{G\"ollnitz-Gordon marking}

The definition of G\"ollnitz-Gordon marking was given in \cite[Definition 3.1]{He-Zhao-2023}.

 \begin{defi}[G\"ollnitz-Gordon  marking]\label{Gollnitz-Gordon-marking}
The G\"ollnitz-Gordon  marking $GG(\pi)$ of a partition $\pi=(\pi_1,\pi_2,\ldots,\pi_\ell)$ is an assignment of positive integers {\rm(}marks{\rm)} to the parts of  $\pi$ from  smallest to  largest such that the marks are as small as possible subject to the condition that for $1\leq i\leq \ell,$
the integer assigned to $\pi_i$ is different from the integers assigned to the parts $\pi_g$ such that $g>i$ and $\pi_i-\pi_g\leq 2$ with strict inequality if $\pi_i$ is odd.
\end{defi}

For example, the G\"ollnitz-Gordon marking of
\begin{equation}\label{gg-r-1}
\pi=(38,38,36,34,32,30,26,26,22,22,22,18,16,16,14,12,12,10,9,6,6,6,2,1)
\end{equation}
 is
\[\begin{split}GG(\pi)=(&38_3,38_1,36_2,34_1,32_2,30_1,26_2,26_1,22_3,22_2,22_1,18_2,16_3,16_1,\\
                        &14_2,12_3,12_1,10_2,9_1,6_3,6_2,6_1,2_2,1_1),
\end{split}\]
where the subscript of each part represents the mark in the G\"ollnitz-Gordon marking.

The G\"ollnitz-Gordon marking of a partition can be represented by an array, where the column indicates the size of a part and the row (counted from bottom to top) indicates the mark, so the G\"ollnitz-Gordon marking of $\pi$ defined in \eqref{gg-r-1} would be
\begin{equation}\label{GG-example-1}
GG(\pi)=\left[
\begin{array}{cccccccccccccccccccc}
&&6&&&12&&16&&22&&&&&&38\\
&2&6&&10&&14&&18&22&26&&32&&36&\\
1&&6&9&&12&&16&&22&26&30&&34&&38
\end{array}
\right].
\end{equation}

For $i\geq 1$, let $N_i(\pi)$ (or $N_i$ for short) denote the number of  parts in the $i$-th row of $GG(\pi)$. From the definition of G\"ollnitz-Gordon marking, it is not hard to find that $N_1\geq N_2\geq\cdots$. We use $\pi^{(i)}=(\pi^{(i)}_1,\pi^{(i)}_2,\ldots,\pi^{(i)}_{N_i})$ to denote the sub-partition of $\pi$ that consists of all $i$-marked parts in $GG(\pi)$. For convention, we define $\pi^{(i)}_0=+\infty$ and $\pi^{(i)}_{N_i+1}=-\infty$.

Let $\pi$ be the partition with G\"ollnitz-Gordon marking given in \eqref{GG-example-1}. By definition, we have $\pi^{(1)}=(38,34,30,26,22,16,12,9,6,1)$, $\pi^{(2)}=(36,32,26,22,18,14,10,6,2)$ and $\pi^{(3)}=(38,22,16,12,6)$. So, we get $N_1(\pi)=10$, $N_2(\pi)=9$ and $N_3(\pi)=5$.

For $k\geq r\geq1$, let $\mathbb{C}(k,r;n)$ denote the set of partitions counted by $B(1;2,k,r;n)$. Define
\[\mathbb{C}(k,r)=\bigcup_{n\geq0}\mathbb{C}(k,r;n).\]
 Clearly, a partition $\pi$  is in  $\mathbb{C}(k,r)$  if and only if no odd part is
repeated,  the marks  of $1$ and $2$  are not exceed to $r-1$  and there are at most $k-1$ rows in $GG(\pi)$.

\subsection{Starting types}

In the rest of this article, we fix $k\geq r\geq 3$. For  $N_2\geq 1$, let $\pi$ be a partition in $\mathbb{C}(k,r)$ such that there are $N_2$ parts marked with $2$ in $GG(\pi)$. We define the starting types of $\pi^{(2)}_1,\pi^{(2)}_{2},\ldots,\pi^{(2)}_{N_2}$ based on the G\"ollnitz-Gordon marking of $\pi$.

{\noindent \bf Starting types:} Assume that $l$ is the largest integer such that there does not exist odd part of $\pi$ greater than or equal to $\pi^{(2)}_l$. Then we say that $\pi^{(2)}_{i}$ is of type $s_{-1}$ for $l+1\leq i\leq N_2$. If $l=0$, then we are done. If $l\geq 1$, then we define the starting types of $\pi^{(2)}_1,\pi^{(2)}_{2},\ldots,\pi^{(2)}_{l}$ from largest to smallest.

We first define the starting type of $\pi^{(2)}_{1}$ as follows.

\begin{itemize}

\item[]Case 1: There is a $1$-marked $\pi^{(2)}_1-1$ in $GG(\pi)$ and  $\pi^{(2)}_1+2$ does not occur in $\pi$. We  say that $\pi^{(2)}_{1}$ is of starting type $s_0$ and we set $s_1(\pi)=\pi^{(2)}_{1}-1$.

\item[]Case 2: There is a $1$-marked $\pi^{(2)}_1-2$ in $GG(\pi)$ and  $\pi^{(2)}_1+2$ does not occur in $\pi$. We say that $\pi^{(2)}_{1}$ is of starting type $s_1$ and we set $s_1(\pi)=\pi^{(2)}_{1}-2$.

\item[]Case 3: There is a $1$-marked $\pi^{(2)}_1+2$ in $GG(\pi)$. We say that $\pi^{(2)}_{1}$ is of starting type $s_2$ and we set $s_1(\pi)=\pi^{(2)}_{1}+2$.

\item[]Case 4: There is a $1$-marked $\pi^{(2)}_1$ in $GG(\pi)$. We say that $\pi^{(2)}_{1}$ is of starting type $s_3$ and we set $s_1(\pi)=\pi^{(2)}_{1}$.

\end{itemize}

For $\pi^{(2)}_{2},\ldots,\pi^{(2)}_{l}$, we set $b=2$ and  repeat the following procedure until $b=l+1$:

(A) We define the starting type of $\pi^{(2)}_{b}$ as follows.

\begin{itemize}
\item[]Case 1:  There is a $1$-marked $\pi^{(2)}_b-1$ in $GG(\pi)$, and $s_{b-1}(\pi)=\pi^{(2)}_b+2$ if there is a $1$-marked $\pi^{(2)}_b+2$ in $GG(\pi)$. We say that $\pi^{(2)}_{b}$ is of starting type $s_0$ and we set $s_b(\pi)=\pi^{(2)}_{b}-1$.

\item[]Case 2:  There is a $1$-marked $\pi^{(2)}_b-2$ in $GG(\pi)$, and $s_{b-1}(\pi)=\pi^{(2)}_b+2$ if there is a $1$-marked $\pi^{(2)}_b+2$ in $GG(\pi)$. We say that $\pi^{(2)}_{b}$ is of starting type $s_1$ and we set $s_b(\pi)=\pi^{(2)}_{b}-2$.

\item[]Case 3:  There is a $1$-marked $\pi^{(2)}_b+2$ in $GG(\pi)$ and $s_{b-1}(\pi)\neq\pi^{(2)}_b+2$. We say that $\pi^{(2)}_{b}$ is of starting type $s_2$ and we set $s_b(\pi)=\pi^{(2)}_{b}+2$.

\item[]Case 4: There is a $1$-marked $\pi^{(2)}_b$ in $GG(\pi)$. We say that $\pi^{(2)}_{b}$ is of starting type $s_3$ and we set $s_b(\pi)=\pi^{(2)}_{b}$.

\end{itemize}

(B) Replace $b$ by $b+1$. If $b=l+1$, then we are done. Otherwise, go back to (A).

For example, let $\pi$ be the partition with G\"ollnitz-Gordon marking given in \eqref{GG-example-1}. It is clear that $\pi\in\mathbb{C}(4,3)$, $N_2=9$ and $l=7$. So, $\pi^{(2)}_{8}=6$ and $\pi^{(2)}_{9}=2$ are of starting type $s_{-1}$. Then, it can be checked that $\pi^{(2)}_{1}=36$ and $\pi^{(2)}_{2}=32$ are of starting type $s_{2}$, $\pi^{(2)}_{3}=26$ and $\pi^{(2)}_{4}=22$ are of starting type $s_{3}$,
$\pi^{(2)}_{5}=18$ and $\pi^{(2)}_{6}=14$ are of starting type $s_{1}$ and $\pi^{(2)}_{7}=10$ is of starting type $s_{0}$.

\subsection{The set $\mathbb{C}_{<}(k,r|p,t)$}

In the remaining of this article, we assume that $p$ and $t$ are integer such that $p,t\geq0$. We give the explicit definition of $\mathbb{C}_{<}(k,r|p,t)$ in Theorem \ref{equiv-main-lem-1-1}.

\begin{defi}\label{defi-c<}

Let $\mathbb{C}_{<}(k,r|p,t)$ be the set of partitions $\pi$ in $\mathbb{C}(k,r)$ such that

\begin{itemize}

\item[{\rm(1)}] there is no odd part of $\pi$ greater than or equal to $2t+1${\rm;}

\item[{\rm(2)}] $\pi^{(2)}_{p+1}<2t+1<\pi^{(2)}_p${\rm;}

\item[{\rm(3)}] if $\pi^{(2)}_{p}=2t+2$, then $\pi^{(2)}_{p}$ is of starting type $s_2$ or $s_3${\rm;}

\item[{\rm(4)}] if $\pi^{(2)}_{p+1}=2t$, then $\pi^{(2)}_{p+1}$ is of starting type $s_0$ or $s_1${\rm.}

\end{itemize}

\end{defi}

For example, let $\pi$ be the partition defined in \eqref{GG-example-1}. It can be checked that $\pi$ is a partition in $\mathbb{C}_{<}(4,3|6,5)$, $\mathbb{C}_{<}(4,3|5,7)$, $\mathbb{C}_{<}(4,3|4,9)$, $\mathbb{C}_{<}(4,3|4,10)$, $\mathbb{C}_{<}(4,3|3,12)$, $\mathbb{C}_{<}(4,3|2,14)$, $\mathbb{C}_{<}(4,3|2,15)$, $\mathbb{C}_{<}(4,3|1,17)$ and $\mathbb{C}_{<}(4,3|0,19)$. But, $\pi$ is not a partition in $\mathbb{C}_{<}(4,3|6,6)$, $\mathbb{C}_{<}(4,3|5,8)$, $\mathbb{C}_{<}(4,3|3,11)$, $\mathbb{C}_{<}(4,3|2,13)$, $\mathbb{C}_{<}(4,3|1,16)$ and $\mathbb{C}_{<}(4,3|0,18)$.

The following proposition is a consequence of the condition (2) in Definition \ref{defi-c<}.
\begin{prop}\label{geqN2}
For $N_2\geq 0$, let $\pi$ be a partition in $\mathbb{C}_{<}(k,r|p,t)$ such that there are $N_2$ parts marked with $2$ in $GG(\pi)$. Then, we have $p+t\geq N_2$.
\end{prop}

\pf If $p=N_2$, then the proposition is obviously right. If $p<N_2$, then by the condition (2) in Definition \ref{defi-c<}, we get
\[2t+1>\pi^{(2)}_{p+1}\geq \pi^{(2)}_{p+2}+2\geq\cdots\geq \pi^{(2)}_{N_2}+2(N_2-p-1)\geq 1+2(N_2-p-1).\]
It yields $t>N_2-p-1$, and so $p+t\geq N_2$. This completes the proof.  \qed

 We will divide  $\mathbb{C}_{<}(k,r|p,t)$ into twelve disjoint subsets and investigate the properties of them. Before doing this, we give the following lemma, which will be related the subsets $\mathbb{C}^{(1)}_{<}(k,r|p,t)$, $\mathbb{C}^{(2)}_{<}(k,r|p,t)$ and $\mathbb{C}^{(3)}_{<}(k,r|p,t)$ of $\mathbb{C}_{<}(k,r|p,t)$.

\begin{lem}\label{lem>2t+6}
Let $\pi$ be a partition in $\mathbb{C}_{<}(k,r|p,t)$ such that $\pi^{(2)}_{p}\geq 2t+6$. Then, the marks of parts $2t+2$ and $2t+4$ are at most $1$  in $GG(\pi)$.
	
\end{lem}

\pf  Suppose to the contrary that there exist parts $2t+2$ or $2t+4$ with marks greater than $1$ in $GG(\pi)$. By the condition (2) in Definition \ref{defi-c<}, we have $\pi^{(2)}_{p+1}<2t+1$. Under the condition that $\pi^{(2)}_{p}\geq 2t+6$, we see that there is no $2$-marked $2t+2$ and $2t+4$ in $GG(\pi)$.

It follows from the definition of G\"ollnitz-Gordon marking that there do not exist parts $2t+4$ with marks greater than $2$ in $GG(\pi)$.  Therefore, there exist parts $2t+2$ with marks greater than $2$ in $GG(\pi)$.
Moreover, there is a $1$-marked $2t$ or $2t+2$ in $GG(\pi)$ and there is a $2$-marked $2t$  in $GG(\pi)$. So, we obtain that $\pi^{(2)}_{p+1}=2t$ and it is of starting type $s_2$ or $s_3$, which contradicts the condition (4) in Definition \ref{defi-c<}. The proof is complete.   \qed

The following corollary immediately follows from Lemma \ref{lem>2t+6}.
\begin{core}
Let $\pi$ be a partition in $\mathbb{C}_{<}(k,r|p,t)$ such that $\pi^{(2)}_{p}\geq 2t+6$. Then,  $2t+2$ and $2t+4$ can not both occur in $\pi$.
\end{core}

For easier expression, we introduce starting cluster indexes based on starting types.

\begin{defi}
For $p\geq 1$, let $\pi$ be a partition in $\mathbb{C}_{<}(k,r|p,t)$. The starting cluster indexes  of $\pi^{(2)}_1,\pi^{(2)}_{2},\ldots,\pi^{(2)}_{p}$ are defined as follows.

Set $b=0$ and $p_0=p+1$, we do the following process.

\begin{itemize}

\item[{\rm(A)}] Assume that $p_{b+1}$ is the smallest integer such that
\[\pi^{(2)}_{p_{b+1}}=\pi^{(2)}_{p_{b}-1}+4(p_b-1-p_{b+1}),\]
and $\pi^{(2)}_{p_{b}-1},\ldots,\pi^{(2)}_{p_{b+1}}$ are of the same starting type. We say that $p_{b+1}$ is the $(b+1)$-th starting cluster index of $\pi$.

\item[{\rm(B)}] Replace $b$ by $b+1$. If $p_b=1$, then we are done. Otherwise, go back to  {\rm(A)}.

\end{itemize}

\end{defi}

For example, let $\pi$ be the partition with G\"ollnitz-Gordon marking given in \eqref{GG-example-1}. If $p=6$ and $t=5$, then we have $p_1=5$, $p_2=3$ and $p_3=1$. If $p=5$ and $t=7$, then we also have $p_1=5$, $p_2=3$ and $p_3=1$. If $p=3$ and $t=12$, then we have $p_1=3$ and $p_2=1$.

We need the following proposition, which will be related the subsets  $\mathbb{C}^{(10)}_{<}(k,r|p,t)$, $\mathbb{C}^{(11)}_{<}(k,r|p,t)$ and $\mathbb{C}^{(12)}_{<}(k,r|p,t)$ of $\mathbb{C}_{<}(k,r|p,t)$.

\begin{prop}\label{prop-i-2}
Let $\pi$ be a partition in  $\mathbb{C}_{<}(k,r|p,t)$ such that $\pi^{(2)}_p=2t+2$ with starting type $s_3$. Then, $\pi^{(2)}_{p_1}+4$ occurs at most once in $\pi$, where $p_1$ is the first starting cluster of $\pi$.
\end{prop}

\pf Suppose to the contrary that $\pi^{(2)}_{p_1}+4$ occurs at least twice  in $\pi$. By the definition of G\"ollnitz-Gordon marking, we see that there are $1$-marked and $2$-marked parts $\pi^{(2)}_{p_1}+4$ in $GG(\pi)$. It yields that $\pi^{(2)}_{p_1-1}=\pi^{(2)}_{p_1}+4$ with starting type $s_3$, which contradicts the choice of $p_1$. So, $\pi^{(2)}_{p_1}+4$ occurs at most once in $\pi$. This completes the proof. \qed

Now, we are in a position to give the twelve disjoint subsets of  $\mathbb{C}_{<}(k,r|p,t)$.

Let $\mathbb{C}^{(1)}_<(k,r|p,t)$ denote the set of partitions $\pi$ in  $\mathbb{C}_{<}(k,r|p,t)$ such that $\pi^{(2)}_{p}\geq 2t+6$, and if $\pi^{(2)}_{p}=2t+6$ then $\pi^{(2)}_{p}$ is of starting type $s_2$ or $s_3$ and the largest mark of parts $2t+6$ in $GG(\pi)$ is $2$.

Let $\mathbb{C}^{(2)}_<(k,r|p,t)$ denote the set of partitions $\pi$ in  $\mathbb{C}_{<}(k,r|p,t)$ such that $\pi^{(2)}_{p}=2t+6$, $2t+2$ does not occur in $\pi$, and if $\pi^{(2)}_{p}$ is of starting type $s_2$ or $s_3$ then there exist parts $2t+6$ with marks greater than $2$ in $GG(\pi)$.

Let $\mathbb{C}^{(3)}_<(k,r|p,t)$ denote the set of partitions $\pi$ in  $\mathbb{C}_{<}(k,r|p,t)$ such that $\pi^{(2)}_{p}=2t+6$ with starting type $s_3$, $2t+2$ occurs in $\pi$ and there exist parts $2t+6$ with marks greater than $2$ in $GG(\pi)$.

Let $\mathbb{C}^{(4)}_<(k,r|p,t)$ denote the set of partitions $\pi$ in  $\mathbb{C}_{<}(k,r|p,t)$ such that $\pi^{(2)}_{p}=2t+4$ with starting type $s_3$ and the largest mark of parts $2t+4$ in $GG(\pi)$ is $2$.

Let $\mathbb{C}^{(5)}_<(k,r|p,t)$ denote the set of partitions $\pi$ in  $\mathbb{C}_{<}(k,r|p,t)$ such that $\pi^{(2)}_{p}=2t+4$ with starting type $s_3$ and there exist parts $2t+4$ with marks greater than $2$ in $GG(\pi)$.

Let $\mathbb{C}^{(6)}_<(k,r|p,t)$ denote the set of partitions $\pi$ in  $\mathbb{C}_{<}(k,r|p,t)$ such that $\pi^{(2)}_{p}=2t+4$ with starting type $s_1$.

Let $\mathbb{C}^{(7)}_<(k,r|p,t)$ denote the set of partitions $\pi$ in  $\mathbb{C}_{<}(k,r|p,t)$ such that $\pi^{(2)}_{p}=2t+4$ with starting type $s_2$.

Let $\mathbb{C}^{(8)}_<(k,r|p,t)$ denote the set of partitions $\pi$ in  $\mathbb{C}_{<}(k,r|p,t)$ such that $\pi^{(2)}_p=2t+2$ with starting type $s_2$.

Let $\mathbb{C}^{(9)}_<(k,r|p,t)$ denote the set of partitions $\pi$ in  $\mathbb{C}_{<}(k,r|p,t)$ such that $\pi^{(2)}_p=2t+2$ with starting type $s_3$ and $\pi^{(2)}_{p_1}+4$ does not occur in $\pi$, where $p_1$ is the first starting cluster index of $\pi$.

Let $\mathbb{C}^{(10)}_<(k,r|p,t)$ denote the set of partitions $\pi$ in  $\mathbb{C}_{<}(k,r|p,t)$ such that $\pi^{(2)}_p=2t+2$ with starting type $s_3$, $\pi^{(2)}_{p_1}+4$ occurs in $\pi$ and $\pi^{(2)}_{p_1}+6$ does not occur in $\pi$, where $p_1$ is the first starting cluster index of $\pi$.

Let $\mathbb{C}_{<}^{(11)}(k,r|p,t)$ denote the set of partitions $\pi$ in  $\mathbb{C}_{<}(k,r|p,t)$ such that $\pi^{(2)}_p=2t+2$ with starting type $s_3$ and  $\pi^{(2)}_{p_1-1}=\pi^{(2)}_{p_1}+6$ with starting type $s_1$, where $p_1$ is the first starting cluster index of $\pi$.

Let $\mathbb{C}_{<}^{(12)}(k,r|p,t)$ denote the set of partitions $\pi$ in  $\mathbb{C}_{<}(k,r|p,t)$ such that $\pi^{(2)}_p=2t+2$ with starting type $s_3$ and  $\pi^{(2)}_{p_1-1}=\pi^{(2)}_{p_1}+6$ with starting type $s_2$, where $p_1$ is the first starting cluster index of $\pi$.

For example, let $\pi$ be the partition  defined in \eqref{GG-example-1}. It can be checked that $\pi$ is a partition in $\mathbb{C}^{(6)}_{<}(4,3|6,5)$, $\mathbb{C}^{(6)}_{<}(4,3|5,7)$, $\mathbb{C}^{(5)}_{<}(4,3|4,9)$, $\mathbb{C}^{(12)}_{<}(4,3|4,10)$, $\mathbb{C}^{(12)}_{<}(4,3|3,12)$, $\mathbb{C}^{(7)}_{<}(4,3|2,14)$, $\mathbb{C}^{(8)}_{<}(4,3|2,15)$, $\mathbb{C}^{(8)}_{<}(4,3|1,17)$ and $\mathbb{C}^{(1)}_{<}(4,3|0,19)$.

Clearly, we have
\[\mathbb{C}_{<}(k,r|p,t)=\bigcup_{1\leq j\leq 12}\mathbb{C}^{(j)}_{<}(k,r|p,t).\]

Then, we give the definition of insertion index.

\begin{defi}
Let $\pi$ be a partition in $\mathbb{C}_{<}(k,r|p,t)$. Assume that $p_1$ and $p_2$ are the first and the second starting cluster indexes of $\pi$ respectively.  We define the insertion index $I_{p,t}(\pi)$ of $\pi$ as follows.

\begin{itemize}
\item[{\rm(1)}] We set $I_{p,t}(\pi)=2t+2$  for $\pi\in\mathbb{C}^{(j)}_<(k,r|p,t)$, where $1\leq j\leq 5$.

 \item[{\rm(2)}] We set $I_{p,t}(\pi)=\pi^{(2)}_{p_1}$ for $\pi\in\mathbb{C}^{(6)}_<(k,r|p,t)$.

 \item[{\rm(3)}]We set $I_{p,t}(\pi)=\pi^{(2)}_{p_1}+2$  for $\pi\in\mathbb{C}^{(j)}_<(k,r|p,t)$, where $7\leq j\leq 9$.

     \item[{\rm(4)}] We set $I_{p,t}(\pi)=\pi^{(2)}_{p_1}+4$ for $\pi\in\mathbb{C}^{(10)}_<(k,r|p,t)$.

      \item[{\rm(5)}] We set $I_{p,t}(\pi)=\pi^{(2)}_{p_2}$ for $\pi\in\mathbb{C}^{(11)}_<(k,r|p,t)$.

      \item[{\rm(6)}] We set $I_{p,t}(\pi)=\pi^{(2)}_{p_2}+2$ for $\pi\in\mathbb{C}^{(12)}_<(k,r|p,t)$.

\end{itemize}

\end{defi}

For example, let $\pi$ be the partition  with G\"ollnitz-Gordon marking presented in \eqref{GG-example-1}. If $p=6$ and $t=5$, then we have
$I_{p,t}(\pi)=\pi^{(2)}_{p_1}=\pi^{(2)}_{5}=18$. If $p=5$ and $t=7$, then we also have $I_{p,t}(\pi)=\pi^{(2)}_{p_1}=\pi^{(2)}_{5}=18$. If $p=3$ and $t=12$, then we have $I_{p,t}(\pi)=\pi^{(2)}_{p_2}+2=\pi^{(2)}_{1}+2=38$.

We conclude this subsection with properties of insertion index.

\begin{prop}\label{prop-i-1}
Let $\pi$ be a partition in $\mathbb{C}_{<}(k,r|p,t)$. Then, $I_{p,t}(\pi)$ and $I_{p,t}(\pi)+2$ can not both occur in $\pi$. More precisely,
\begin{itemize}
\item[{\rm(1)}]  the marks of parts $2t+2$ and $2t+4$ are at most $1$ in $GG(\pi)$ for $\pi\in\mathbb{C}^{(j)}_<(k,r|p,t)$, where $1\leq j\leq 3${\rm;}

 \item[{\rm(2)}] $2t+2$ does not occur in $\pi$ for $\pi\in\mathbb{C}^{(j)}_<(k,r|p,t)$, where $j=4,5${\rm;}

 \item[{\rm(3)}] $I_{p,t}(\pi)+2$ does not occur in $\pi$ for $\pi\in\mathbb{C}^{(j)}_<(k,r|p,t)$, where $6\leq j\leq 12$.

\end{itemize}

\end{prop}

\pf (1) It is an immediate consequence of Lemma \ref{lem>2t+6}.

(2) Suppose to the contrary that $2t+2$ occurs in $\pi$. By the definitions of $\mathbb{C}^{(4)}_<(k,r|p,t)$ and $\mathbb{C}^{(5)}_<(k,r|p,t)$, we know that $\pi^{(2)}_{p}=2t+4$ and it is of starting type $s_3$, and so there are $1$-marked and $2$-marked parts $2t+4$ in $GG(\pi)$.
By virtue of the definition of G\"ollnitz-Gordon marking, we find that  the marks of parts $2t+2$ are greater than $2$ and there are $1$-marked and $2$-marked parts $2t$ in $GG(\pi)$. It implies that  $\pi^{(2)}_{p+1}=2t$ and it is of starting type $s_3$, which contradicts the condition (4) in Definition \ref{defi-c<}. Hence, $2t+2$ does not occur in $\pi$.

(3) Appealing to the definitions of starting type and starting cluster index, we see that $I_{p,t}(\pi)+2$ does not occur in $\pi$ for $\pi\in\mathbb{C}^{(j)}_<(k,r|p,t)$, where $j=6,7,8,11,12$.

If $\pi\in\mathbb{C}^{(9)}_<(k,r|p,t)$, then we have $I_{p,t}(\pi)=\pi^{(2)}_{p_1}+2$ and $\pi^{(2)}_{p_1}+4$ does not occur in $\pi$. It implies that  $I_{p,t}(\pi)+2$ does not occur in $\pi$.

If $\pi\in\mathbb{C}^{(10)}_<(k,r|p,t)$, then we have $I_{p,t}(\pi)=\pi^{(2)}_{p_1}+4$ and $\pi^{(2)}_{p_1}+6$ does not  occur in $\pi$. It yields that $I_{p,t}(\pi)+2$ does not occur in $\pi$.

We can conclude that the condition (3) is satisfied. Thus, we complete the proof.   \qed

\subsection{The set $\mathbb{C}_{=}(k,r|p,t)$}

In this subsection, we will give the explicit definition of    $\mathbb{C}_{=}(k,r|p,t)$ in Theorem \ref{equiv-main-lem-1-1}. Before doing this, we need the following proposition.

\begin{lem}\label{prop-new-add}
Let $\pi$ be a partition in $\mathbb{C}(k,r)$ such that $2t+1$ and $2t+2$ both occur in $\pi$. Then, the marks of parts $2t+2$ in $GG(\pi)$ are greater than the mark of $2t+1$ in   $GG(\pi)$.
\end{lem}	

\pf Assume that $2t+1$ is marked with $r$ in $GG(\pi)$. It follows from the definition of G\"ollnitz-Gordon marking that $2t+2$ can not be marked with $r$ in $GG(\pi)$. We consider the following two cases.

Case 1: If $r=1$, then it is obviously right.

Case 2: If $r\geq 2$, then by the definition of G\"ollnitz-Gordon marking, we see that there are $1$-marked, $2$-marked, \ldots, $(r-1)$-marked parts $2t$ in $GG(\pi)$. It implies that  there are no $1$-marked, $2$-marked, \ldots, $(r-1)$-marked parts $2t+2$ in $GG(\pi)$. So, the marks of parts $2t+2$ in $GG(\pi)$ are greater than $r$.

Thus, we have completed the proof.   \qed

The following corollary is an immediate consequence of Lemma \ref{prop-new-add}.
\begin{core}\label{2t+2}
Let $\pi$ be a partition in $\mathbb{C}(k,r)$ with the largest odd part $2t+1$. Then,
\begin{itemize}
\item[\rm{(1)}] if there is a $2$-marked $2t+2$ in $GG(\pi)$, then it is of starting type $s_0$ or $s_2${\rm;}

\item[\rm{(2)}] if there is a $2$-marked $2t+4$ in $GG(\pi)$, then it is of starting type $s_3$.

\end{itemize}
\end{core}

\pf In view of Lemma \ref{prop-new-add}, we find that there is no $1$-marked $2t+2$ in $GG(\pi)$, and if there is a $2$-marked $2t+2$ in $GG(\pi)$ then $2t+1$ is marked with $1$ in $GG(\pi)$.
 It follows from the definitions of G\"ollnitz-Gordon marking and starting type that the conditions (1) and (2) are verified. The proof is complete.  \qed

Now, we proceed to introduce the definition of $\mathbb{C}_{=}(k,r|p,t)$.

\begin{defi}\label{defi-=}

Let $\mathbb{C}_{=}(k,r|p,t)$ be the set of partitions $\pi$ in $\mathbb{C}(k,r)$ such that

\begin{itemize}

\item[{\rm(1)}] the largest odd part of $\pi$ is $2t+1${\rm;}

\item[{\rm(2)}] the mark of $2t+1$  in $GG(\pi)$ is at most $2${\rm;}

\item[{\rm(3)}] $\pi^{(2)}_{p}\geq2t+2$ and $\pi^{(2)}_{p+1}\leq2t+2${\rm;}

\item[{\rm(4)}] if there is a $2$-marked $2t+2$ in $GG(\pi)$ and it is of starting type $s_0$, then $\pi^{(2)}_{p+1}=2t+2$ and there exists $i$ such that $i\leq p+1$, $\pi^{(2)}_{i}=\pi^{(2)}_{p+1}+4(p-i+1)$ and $\pi^{(2)}_{i}$ occurs once in $\pi${\rm;}

\item[{\rm(5)}] if there is a $2$-marked $2t+2$ in $GG(\pi)$ and it is of starting type $s_2$, then $\pi^{(2)}_{p}=2t+2${\rm;}

\item[{\rm(6)}] if $2t+2$ occurs in $\pi$ and there is no $2$-marked $2t+2$ in $GG(\pi)$, then $\pi^{(2)}_{p}=2t+4$ with starting type $s_3$ and there exists $i$ such that $i\leq p$,  $\pi^{(2)}_{i}=\pi^{(2)}_{p}+4(p-i)$ and $\pi^{(2)}_{i}+2$ does not occur in $\pi$.

\end{itemize}

\end{defi}

For example, let $\pi$ be the partition in $\mathbb{C}(4,3)$, whose G\"ollnitz-Gordon marking is given in \eqref{GG-example-1}. The largest odd part of $\pi$ is $9$, which is marked with $1$ in $GG(\pi)$. We find that $\pi^{(2)}_{7}=10$ is of starting type $s_0$ and $\pi^{(2)}_{7}=10$ occurs once in $\pi$. It yields $\pi\in\mathbb{C}_{=}(4,3|6,4)$.

For another example, let $\pi$ be the partition in $\mathbb{C}(4,3)$ with G\"ollnitz-Gordon marking
\begin{equation}\label{GG-example-newlabel-10}
GG(\pi)=\left[
\begin{array}{cccccccccccccccccccc}
&&6&&&&12&&16&&&24&&&&&38\\
&2&6&&10&&&14&&18&22&&&28&&34&38\\
1&&6&9&&11&&14&&18&22&&26&&30&34&38
\end{array}
\right].
\end{equation}
The largest odd part of $\pi$ is $11$, which is marked with $1$ in $GG(\pi)$. We see that $12$ occurs in $\pi$, there is no $2$-marked $12$ in $GG(\pi)$ and $\pi^{(2)}_{6}=14$. Moreover, it can be checked that $\pi^{(2)}_{5}=18=\pi^{(2)}_{6}+4$, and $\pi^{(2)}_{5}+2=20$ does not occur in $\pi$. Then, we have $\pi\in\mathbb{C}_{=}(4,3|6,5)$.

We proceed to divide  $\mathbb{C}_{=}(k,r|p,t)$ into twelve disjoint subsets and investigate the properties of them. Before doing this, we give the following lemma, which will be related the subsets $\mathbb{C}^{(3)}_{=}(k,r|p,t)$, $\mathbb{C}^{(7)}_{=}(k,r|p,t)$, $\mathbb{C}^{(10)}_{=}(k,r|p,t)$ and $\mathbb{C}^{(12)}_{=}(k,r|p,t)$ of $\mathbb{C}_{=}(k,r|p,t)$.

\begin{lem}\label{rpt-lem-1}
Assume that $\pi$ is a partition in $\mathbb{C}_{=}(k,r|p,t)$, $\pi^{(2)}_{p+1}=2t+2$, and $s$ is the smallest integer such that  $\pi^{(2)}_{s}=\pi^{(2)}_{p+1}+4(p-s+1)$ with starting type $s_0$ or $s_1$.  Then,

\begin{itemize}
	\item[\rm{(1)}]  $\pi^{(2)}_{s}+2$ does not occur in $\pi${\rm;}
	
	\item[\rm{(2)}]  for $i<s$, $\pi^{(2)}_{i}$ is of starting type $s_3$ if  $\pi^{(2)}_{i}=\pi^{(2)}_{p+1}+4(p-i+1)$.
\end{itemize}
	
\end{lem}

\pf It follows from the definition of $\mathbb{C}_{=}(k,r|p,t)$ that  $\pi^{(2)}_{p+1}=2t+2$ is of starting type $s_0$. By the choice of $s$,  we have $s\leq p+1$.

(1) Suppose to the contrary that  $\pi^{(2)}_{s}+2$ occurs in $\pi$. Under the condition that $\pi^{(2)}_{s}$ is of starting type $s_0$ or $s_1$, we see that there is no $1$-marked $\pi^{(2)}_{s}$ in $GG(\pi)$, and so there is a $1$-marked  $\pi^{(2)}_{s}+2$ in $GG(\pi)$. Moreover,  we have $\pi^{(2)}_{s-1}=\pi^{(2)}_{s}+4$ and it is of starting type $s_1$, which contradicts the choice of $s$.
Hence, $\pi^{(2)}_{s}+2$ does not occur in $\pi$.

(2) We just need to show that if  $\pi^{(2)}_{s-1}=\pi^{(2)}_{s}+4$ then  $\pi^{(2)}_{s-1}$  is of starting type $s_3$. Assume that  $\pi^{(2)}_{s-1}=\pi^{(2)}_{s}+4$.
Using the condition (1), we deduce that  there is no $1$-marked $\pi^{(2)}_{s}+2$ in $GG(\pi)$. By the definition of G\"ollnitz-Gordon marking, we see that there is a $1$-marked part $\pi^{(2)}_{s}+4$ in  $GG(\pi)$. It yields that  $\pi^{(2)}_{s-1}$ is of starting type $s_3$. The proof is complete. \qed

Next, we divide  $\mathbb{C}_{=}(k,r|p,t)$ into the following twelve disjoint subsets.

Let $\mathbb{C}_{=}^{(1)}(k,r|p,t)$ denote the set of partitions $\pi$ in  $\mathbb{C}_{=}(k,r|p,t)$ such that $\pi^{(2)}_{p}\geq 2t+8$.

Let $\mathbb{C}_{=}^{(2)}(k,r|p,t)$ denote the set of partitions $\pi$ in  $\mathbb{C}_{=}(k,r|p,t)$ such that $\pi^{(2)}_{p}=2t+6$ with starting type $s_2$ or $s_3$, $\pi^{(2)}_{p+1}<2t+2$, and if $\pi^{(2)}_{p}$ is of starting type $s_2$ and $\pi^{(2)}_{i}=\pi^{(2)}_{p}+4(p-i)$ then $\pi^{(2)}_{i}+2$ occurs at least twice in $\pi$ for $i\leq p$.

Let $\mathbb{C}_{=}^{(3)}(k,r|p,t)$ denote the set of partitions $\pi$ in  $\mathbb{C}_{=}(k,r|p,t)$ such that $\pi^{(2)}_{p}=2t+6$ with starting type $s_3$, $\pi^{(2)}_{p+1}=2t+2$, and if
$\pi^{(2)}_i=\pi^{(2)}_{p}+4(p-i)$ then  $\pi^{(2)}_i+2$ occurs in $\pi$ for $i\leq p$.

Let $\mathbb{C}_{=}^{(4)}(k,r|p,t)$ denote the set of partitions $\pi$ in  $\mathbb{C}_{=}(k,r|p,t)$ such that $\pi^{(2)}_{p}=2t+6$ with starting type $s_1$ or $s_2$, $\pi^{(2)}_{p+1}<2t+2$, and if $\pi^{(2)}_{p}$ is of starting type $s_2$ then there exists $i$ such that $i\leq p$, $\pi^{(2)}_{i}=\pi^{(2)}_{p}+4(p-i)$ and  $\pi^{(2)}_{i}+2$ occurs once in $\pi$.

Let $\mathbb{C}_{=}^{(5)}(k,r|p,t)$ denote the partitions $\pi$ in  $\mathbb{C}_{=}(k,r|p,t)$ such that $\pi^{(2)}_{p}=2t+4$ with starting type $s_3$, and $\pi^{(2)}_{i}+2$ occurs in $\pi$ if $\pi^{(2)}_{i}=\pi^{(2)}_{p}+4(p-i)$ for $i\leq p$.

Let $\mathbb{C}_{=}^{(6)}(k,r|p,t)$ denote the set of partitions $\pi$ in  $\mathbb{C}_{=}(k,r|p,t)$ such that $\pi^{(2)}_{p}=2t+4$ with starting type $s_3$, $2t+1$ is marked with $1$ in $GG(\pi)$,
 and there exists $i$ such that $i\leq p$, $\pi^{(2)}_{i}=\pi^{(2)}_{p}+4(p-i)$ and $\pi^{(2)}_{i}+2$ does not occur in $\pi$.

Let $\mathbb{C}_{=}^{(7)}(k,r|p,t)$ denote the set of partitions $\pi$ in  $\mathbb{C}_{=}(k,r|p,t)$ such that $\pi^{(2)}_{p}=2t+6$, $\pi^{(2)}_{p+1}=2t+2$, $s=p+1$ is the smallest integer such that
$\pi^{(2)}_s=\pi^{(2)}_{p+1}+4(p-s+1)$ and  $\pi^{(2)}_s$ occurs once in $\pi$, and
there  exists $i$ such that $i<p+1$, $\pi^{(2)}_i=\pi^{(2)}_{p+1}+4(p-i+1)$ and  $\pi^{(2)}_i+2$ does not occur in $\pi$.

Let $\mathbb{C}_{=}^{(8)}(k,r|p,t)$ denote the set of partitions $\pi$ in  $\mathbb{C}_{=}(k,r|p,t)$ such that $\pi^{(2)}_{p}=2t+4$  with starting type $s_3$, $2t+1$ is marked with $2$ in $GG(\pi)$,
and there exists $i$ such that $i\leq p$, $\pi^{(2)}_{i}=\pi^{(2)}_{p}+4(p-i)$ and $\pi^{(2)}_{i}+2$ does not occur in $\pi$.

Let $\mathbb{C}_{=}^{(9)}(k,r|p,t)$ denote the set of partitions $\pi$ in  $\mathbb{C}_{=}(k,r|p,t)$ such that $\pi^{(2)}_p=2t+2$, and if
$\pi^{(2)}_i=\pi^{(2)}_{p}+4(p-i)+2$ then  $\pi^{(2)}_i$ is of starting type $s_3$ and $\pi^{(2)}_i+2$ occurs in $\pi$ for $i<s$, where
$s$ is the smallest integer such that $\pi^{(2)}_s=\pi^{(2)}_{p}+4(p-s)$.

Let $\mathbb{C}_{=}^{(10)}(k,r|p,t)$ denote the set of partitions $\pi$ in  $\mathbb{C}_{=}(k,r|p,t)$ such that $\pi^{(2)}_{p}=2t+6$, $\pi^{(2)}_{p+1}=2t+2$, and if
$\pi^{(2)}_i=\pi^{(2)}_{p}+4(p-i)$ then $\pi^{(2)}_i$ is of starting type $s_3$ and $\pi^{(2)}_i+2$ occurs in $\pi$ for $i<s$, where
$s(\leq p)$ is the smallest integer such that $\pi^{(2)}_s=\pi^{(2)}_{p}+4(p-s)$ and  $\pi^{(2)}_s$ occurs once in $\pi$.

Let $\mathbb{C}_{=}^{(11)}(k,r|p,t)$ denote the set of partitions $\pi$ in  $\mathbb{C}_{=}(k,r|p,t)$ such that $\pi^{(2)}_p=2t+2$, and there exists $i$ such that $i<s$, $\pi^{(2)}_i=\pi^{(2)}_{p}+4(p-i)+2$, $\pi^{(2)}_i$ is of starting type $s_3$ and  $\pi^{(2)}_i+2$ does not occur in $\pi$, where
$s$ is the smallest integer such that $\pi^{(2)}_s=\pi^{(2)}_{p}+4(p-s)$.

Let $\mathbb{C}_{=}^{(12)}(k,r|p,t)$ denote the set of partitions $\pi$ in  $\mathbb{C}_{=}(k,r|p,t)$ such that $\pi^{(2)}_{p}=2t+6$, $\pi^{(2)}_{p+1}=2t+2$, and there exists $i$ such that $i<s$, $\pi^{(2)}_i=\pi^{(2)}_{p}+4(p-i)$ and  $\pi^{(2)}_i+2$ does not occur in $\pi$, where
$s(\leq p)$ is the smallest integer such that $\pi^{(2)}_s=\pi^{(2)}_{p}+4(p-s)$ and  $\pi^{(2)}_s$ occurs once in $\pi$.

For example, let $\pi$ be the partition in $\mathbb{C}_{=}(4,3|6,4)$ with G\"ollnitz-Gordon marking stated in \eqref{GG-example-1}. Then, we have $\pi\in\mathbb{C}^{(12)}_{=}(4,3|6,4)$.

For another example, let $\pi$ be the partition in $\mathbb{C}_{=}(4,3|6,5)$, whose G\"ollnitz-Gordon marking is in \eqref{GG-example-newlabel-10}. Then, we have $\pi\in\mathbb{C}^{(6)}_{=}(4,3|6,5)$.

Clearly, we have
\[\mathbb{C}_{=}(k,r|p,t)=\bigcup_{1\leq j\leq 12}\mathbb{C}^{(j)}_{=}(k,r|p,t).\]

Then, we give the definition of division index.

\begin{defi}
Let $\pi$ be a partition in $\mathbb{C}_{=}(k,r|p,t)$.  We define the division index $D_{p,t}(\pi)$ of $\pi$ as follows.
	
	\begin{itemize}
		\item[{\rm(1)}] We set $D_{p,t}(\pi)=2t+2$ for $\pi\in\mathbb{C}^{(j)}_=(k,r|p,t)$, where $1\leq j\leq 5$.

		\item[{\rm(2)}] We set $D_{p,t}(\pi)=\pi^{(2)}_s$, where $s$ is the smallest integer such that $s\leq p$, $\pi^{(2)}_{s}=\pi^{(2)}_{p}+4(p-s)$, and satisfies one of the following conditions{\rm:}

		\begin{itemize}
			
			\item[{\rm(2.1)}]  $\pi^{(2)}_{s}+2$ does not occur in $\pi$ for $\pi\in\mathbb{C}^{(j)}_=(k,r|p,t)$, where $j=6,7,8,12${\rm;}

			\item[{\rm(2.2)}] $\pi^{(2)}_s$ occurs once in $\pi$  for $\pi\in\mathbb{C}^{(10)}_=(k,r|p,t)$.
			
\end{itemize}
		
		\item[{\rm(3)}] We set $D_{p,t}(\pi)=\pi^{(2)}_{s}+2$, where $s$ is the smallest integer such that $s\leq p$ and $\pi^{(2)}_{s}=\pi^{(2)}_{p}+4(p-s)$  for  $\pi\in \mathbb{C}^{(9)}_=(k,r|p,t)$.

\item[{\rm(4)}] We set $D_{p,t}(\pi)=\pi^{(2)}_{s}$, where $s$ is the smallest integer such that $s\leq p$, $\pi^{(2)}_{s}=\pi^{(2)}_{p}+4(p-s)+2$ and $\pi^{(2)}_{s}+2$ does not occur in $\pi$  for  $\pi\in \mathbb{C}^{(11)}_=(k,r|p,t)$.
		
	\end{itemize}
	
\end{defi}

For example, let $\pi$ be the partition in $\mathbb{C}^{(12)}_{=}(4,3|6,4)$ defined in \eqref{GG-example-1}. Then, we have $D_{6,4}(\pi)=\pi^{(2)}_3=26$.

For another example, let $\pi$ be the partition in $\mathbb{C}^{(6)}_{=}(4,3|6,5)$ with G\"ollnitz-Gordon marking given in \eqref{GG-example-newlabel-10}. Then, we have $D_{6,5}(\pi)=\pi^{(2)}_5=18$.

We conclude this subsection with properties of division index.

\begin{prop}\label{prop-divisio-index-1}
Let $\pi$ be a partition in $\mathbb{C}_{=}(k,r|p,t)$. Then, $D_{p,t}(\pi)$ and $D_{p,t}(\pi)+2$ can not both occur in $\pi$. More precisely,
\begin{itemize}
\item[{\rm(1)}] $2t+2$ can only be marked with $2$ in $GG(\pi)$, $2t+4$ can only be marked with $1$ in $GG(\pi)$, and $2t+2$ and $2t+4$ can not both occur in $\pi$ for $\pi\in\mathbb{C}^{(j)}_=(k,r|p,t)$, where $1\leq j\leq 3${\rm;}

 \item[{\rm(2)}] $2t+2$ does not occur in $\pi$ for $\pi\in\mathbb{C}^{(j)}_=(k,r|p,t)$, where $j=4,5${\rm;}

 \item[{\rm(3)}] $D_{p,t}(\pi)+2$ does not occur in $\pi$ for $\pi\in\mathbb{C}^{(j)}_=(k,r|p,t)$, where $6\leq j\leq 12$.

\end{itemize}

\end{prop}

\pf (1) Let $\pi$ be a partition in $\mathbb{C}^{(1)}_=(k,r|p,t)$ or $\mathbb{C}^{(2)}_=(k,r|p,t)$ or $\mathbb{C}^{(3)}_=(k,r|p,t)$. Then, we have $\pi^{(2)}_p\geq 2t+6$. We proceed to show that

\begin{itemize}
\item[(A)] $2t+2$ can only be marked with $2$ in $GG(\pi)$;

 \item[(B)] $2t+2$ and $2t+4$ can not both occur in $\pi$;

\item[(C)]  $2t+4$ can only be marked with $1$ in $GG(\pi)$.

\end{itemize}

{\noindent Condition (A).} Assume that $2t+2$ occurs in $\pi$. Note that $\pi^{(2)}_p\geq 2t+6$, then by the condition (6) in Definition \ref{defi-=}, we see that there is a $2$-marked $2t+2$ in  $GG(\pi)$. Under the condition (3) in Definition \ref{defi-=}, we know that $\pi^{(2)}_{p+1}\leq 2t+2$, which implies that $\pi^{(2)}_{p+1}=2t+2$. It yields that $\pi\not\in\mathbb{C}^{(2)}_=(k,r|p,t)$, and so we have  $\pi\in\mathbb{C}^{(1)}_=(k,r|p,t)$ or $\mathbb{C}^{(3)}_=(k,r|p,t)$.
 Moreover, we find that if there exists $i$ such that $i<p+1$ and $\pi^{(2)}_i=\pi^{(2)}_{p+1}+4(p-i+1)$ then we have $\pi\in\mathbb{C}^{(3)}_=(k,r|p,t)$ and $\pi^{(2)}_i$ is of starting type $s_3$.
Combining with the condition (4) in Definition \ref{defi-=}, we find that $\pi^{(2)}_{p+1}=2t+2$ occurs once in $\pi$. So, the condition (A) is satisfied.

{\noindent Condition (B).} Suppose to contrary that  $2t+2$ and $2t+4$ both occur in $\pi$. Using the argument in the proof of the condition (A), we know that $\pi^{(2)}_{p+1}=2t+2$. By the condition (4) in Definition \ref{defi-=}, we deduce that $\pi^{(2)}_{p+1}=2t+2$  is of starting type $s_0$. It follows from the definition of G\"ollnitz-Gordon marking that there is a $1$-marked $2t+4$  in $GG(\pi)$. By the definition of starting type, we obtain that $\pi^{(2)}_{p}=2t+6$ and it is of starting type $s_1$, which  leads to a contradiction. So,  the condition (B) is verified.

{\noindent Condition (C).} Assume that there are $r(\geq 1)$ parts $2t+4$ in $\pi$. Using the condition (B), we know that $2t+2$ does not occur in $\pi$. By the definition of G\"ollnitz-Gordon marking, we see that there are $1$-marked, $\ldots$, $r$-marked parts $2t+4$ in  $GG(\pi)$. Under the condition that $\pi^{(2)}_p\geq 2t+6$, we obtain that there is no $2$-marked $2t+4$ in $GG(\pi)$. It implies that $r=1$. This completes the proof of the condition (C).

(2) Let $\pi$ be a partition in $\mathbb{C}^{(4)}_=(k,r|p,t)$ or $\mathbb{C}^{(5)}_=(k,r|p,t)$. Then, there is no $2$-marked $2t+2$ in $GG(\pi)$. Appealing to the condition (6) in Definition \ref{defi-=}, we obtain that   $2t+2$ does not occur in $\pi$.

(3) It is from  the definition of $D_{p,t}(\pi)$ that $D_{p,t}(\pi)+2$ does not occur in $\pi$ for $\pi\in\mathbb{C}^{(j)}_=(k,r|p,t)$, where $j=6,7,8,11,12$.

If $\pi\in\mathbb{C}^{(9)}_=(k,r|p,t)$, then we have $D_{p,t}(\pi)=\pi^{(2)}_{s}+2$, where $s$ is the smallest integer such that $s\leq p$ and $\pi^{(2)}_{s}=\pi^{(2)}_{p}+4(p-s)$. Moreover, we find that
$\pi^{(2)}_{s}$ is of starting type $s_2$. So, we see that $D_{p,t}(\pi)+2$ does not occur in $\pi$, otherwise we obtain that $\pi^{(2)}_{s-1}=D_{p,t}(\pi)+2=\pi^{(2)}_{p}+4(p-s+1)$, which contradicts the choice of $s$.

If $\pi\in\mathbb{C}^{(10)}_=(k,r|p,t)$, then we have $D_{p,t}(\pi)=\pi^{(2)}_s$, where $s$ is the smallest integer such that $s\leq p$, $\pi^{(2)}_{s}=\pi^{(2)}_{p}+4(p-s)$ and $\pi^{(2)}_s$ occurs once in $\pi$. Moreover, we find that
$\pi^{(2)}_{s}$ is of starting type $s_1$. It follows from the definition of $\mathbb{C}^{(10)}_=(k,r|p,t)$ that $s$ is the smallest integer such that  $\pi^{(2)}_{s}=\pi^{(2)}_{p+1}+4(p-s+1)$ with starting type $s_0$ or $s_1$. In view of the condition (1) in Lemma \ref{rpt-lem-1}, we get that $\pi^{(2)}_{s}+2$ does not occur in $\pi$, and so $D_{p,t}(\pi)+2$ does not occur in $\pi$.

 In conclusion, the condition (3) is verified. Thus, we complete the proof.  \qed

 \subsection{Equivalent statements of Theorem \ref{equiv-main-lem-1-1}}

 Now, we turn to give equivalent statements of Theorem \ref{equiv-main-lem-1-1}.
We first introduce the subset $\mathbb{C}_\sim(k,r|p,t)$ of $\mathbb{C}_{<}(k,r|p,t)$. To do this, we need to introduce the definition of reduction types, which is a modification of that given by He and Zhao \cite{He-Zhao-2023}.

Let $\pi$ be a partition in $\mathbb{C}_{<}(k,r|p,t)$. Assume that $l$ is the largest integer such that $\pi^{(2)}_l>I_{p,t}(\pi)$. If $l\geq 1$, then set $b=0$ and $p_0=0$ and carry out the following procedure:
\begin{itemize}

    \item[(A)]  Assume that $p_{b+1}$ is the largest integer such that $p_{b+1}\geq p_b+1$, $\pi^{(2)}_{p_b+1}-\pi^{(2)}_{p_{b+1}}=4(p_{b+1}-1-p_{b})$, and satisfying one of the following conditions:
 \begin{itemize}

\item[(1)] $\pi^{(2)}_{i}$ is of starting type $s_3$ and $\pi^{(2)}_{i}+2$ occurs in  $\pi$ for $p_b+1\leq i\leq p_{b+1}$. We say that $\pi^{(2)}_{p_b+1},\ldots,\pi^{(2)}_{p_{b+1}}$ are of insertion type $\hat{A}_1${\rm;}

\item[(2)] $\pi^{(2)}_{i}$ is of starting type $s_3$ for $p_b+1\leq i\leq p_{b+1}$, $\pi^{(2)}_{p_{b}+1}+2$ does not occurs in $\pi$, and there is no $1$-marked  $\pi^{(2)}_{p_{b+1}}-4$ in $GG(\pi)$. We say that $\pi^{(2)}_{p_b+1},\ldots,\pi^{(2)}_{p_{b+1}}$ are of insertion type $\hat{A}_2${\rm;}

 \item[(3)] $\pi^{(2)}_{i}$ is of starting type $s_3$ for $p_b+1\leq i\leq p_{b+1}$, $\pi^{(2)}_{p_{b}+1}+2$ does not occurs in $\pi$, and there is a $1$-marked  $\pi^{(2)}_{p_{b+1}}-4$ in $GG(\pi)$. We say that $\pi^{(2)}_{p_b+1},\ldots,\pi^{(2)}_{p_{b+1}}$ are of insertion type $\hat{A}_3${\rm;}

\item[(4)] $\pi^{(2)}_{i}$ is of starting type $s_2$ and $\pi^{(2)}_{i}+2$ occurs at least twice in  $\pi$ for $p_b+1\leq i\leq p_{b+1}$. We say that $\pi^{(2)}_{p_b+1},\ldots,\pi^{(2)}_{p_{b+1}}$ are of insertion type $\hat{B}${\rm;}

       \item[(5)]  $\pi^{(2)}_{i}$ is of starting type $s_1$ or $s_2$ for $p_b+1\leq i\leq p_{b+1}$, and $\pi^{(2)}_{p_{b}+1}+2$ occurs at most once in  $\pi$. We say that $\pi^{(2)}_{p_b+1},\ldots,\pi^{(2)}_{p_{b+1}}$ are of insertion type $\hat{C}$.
\end{itemize}

Define $\hat{\alpha}_{p_b+1}(\pi)=\cdots=\hat{\alpha}_{p_{b+1}}(\pi)=\{p_b+1,\ldots,p_{b+1}\}$.

\item[(B)] Replace $b$ by $b+1$. If $p_b=l$, then we are done. If $p_b<l$, then go back to (A).

\end{itemize}

For example, let $\pi$ be a partition in $\mathbb{C}^{(2)}_{<}(4,3|9,0)$ with G\"ollnitz-Gordon marking
\begin{equation}\label{example-sim-1}
{GG}(\pi)=\left[
\begin{array}{cccccccccccccccccccc}
&6&&&12&&&&20&&26&&&34&&40\\
&6&&10&&14&&18&&24&&28&32&&38&&42\\
4&&8&&12&&16&&20&24&&28&32&&38&&42
\end{array}
\right].
\end{equation}
It can be checked that $l=9$ is the largest integer such that $\pi^{(2)}_{l}=6>I_{p,t}(\pi)=2$. Moreover, we get that $\pi^{(2)}_{1}=42$ and $\pi^{(2)}_{2}=38$ are of reduction type $\hat{A}_2$ and $\hat{\alpha}_{1}(\pi)=\hat{\alpha}_{2}(\pi)=\{1,2\}$; $\pi^{(2)}_{3}=32$  is of reduction type $\hat{A}_1$ and $\hat{\alpha}_{3}(\pi)=\{3\}$; $\pi^{(2)}_{4}=28$ and $\pi^{(2)}_{5}=24$ are of  reduction type $\hat{A}_3$ and $\hat{\alpha}_{4}(\pi)=\hat{\alpha}_{5}(\pi)=\{4,5\}$; $\pi^{(2)}_{6}=18$  is of  reduction type $\hat{B}$ and $\hat{\alpha}_{6}(\pi)=\{6\}$; $\pi^{(2)}_{7}=14$, $\pi^{(2)}_{8}=10$ and $\pi^{(2)}_{9}=6$ are of reduction type $\hat{C}$ and $\hat{\alpha}_{7}(\pi)=\hat{\alpha}_{8}(\pi)=\hat{\alpha}_{9}(\pi)=\{7,8,9\}$.

Then, we proceed to give the definition of $\mathbb{C}_\sim(k,r|p,t)$, which is the union of the following twelve subsets of  $\mathbb{C}_{<}(k,r|p,t)$.

Let $\mathbb{C}^{(1)}_\sim(k,r|p,t)$ denote the set of partitions $\pi$ in  $\mathbb{C}_{<}(k,r|p,t)$ such that $\pi^{(2)}_{p}\geq 2t+8$.

Let $\mathbb{C}^{(2)}_\sim(k,r|p,t)$ denote the set of partitions $\pi$ in  $\mathbb{C}_{<}(k,r|p,t)$ such that $\pi^{(2)}_{p}=2t+6$ with reduction type $\hat{A}_1$ or $\hat{A}_2$ or $\hat{B}$, and $2t+2$ does not occur in $\pi$.

Let $\mathbb{C}^{(3)}_\sim(k,r|p,t)$ denote the set of partitions $\pi$ in  $\mathbb{C}_{<}(k,r|p,t)$ such that $\pi^{(2)}_{p}=2t+6$ with reduction type $\hat{A}_1$ and $2t+2$ occurs in $\pi$.

Let $\mathbb{C}^{(4)}_\sim(k,r|p,t)$ denote the set of partitions $\pi$ in  $\mathbb{C}_{<}(k,r|p,t)$ such that $\pi^{(2)}_{p}=2t+6$ with reduction type $\hat{C}$, and $2t+2$ does not occur in $\pi$.

Let $\mathbb{C}^{(5)}_\sim(k,r|p,t)$ denote the set of partitions $\pi$ in  $\mathbb{C}_{<}(k,r|p,t)$ such that $\pi^{(2)}_{p}=2t+4$ with reduction type $\hat{A}_1$.

For $6\leq j\leq 12$, let $\mathbb{C}^{(j)}_\sim(k,r|p,t)$ denote the set of partitions $\pi$ in  $\mathbb{C}^{(j)}_{<}(k,r|p,t)$ such that if $\pi^{(2)}_l=I_{p,t}(\pi)+4$ then
 $\pi^{(2)}_l$ is of reduction type $\hat{A}_1$, where $l$ is the largest integer such that $\pi^{(2)}_l>I_{p,t}(\pi)$.

 For example, let $\pi$ be the partition defined in \eqref{example-sim-1}. Then, we see that $\pi$ is a partition in $\mathbb{C}^{(4)}_{\sim}(4,3|9,0)$.

Define
\[\mathbb{C}_{\sim}(k,r|p,t)=\bigcup_{1\leq j\leq 12}\mathbb{C}^{(j)}_{\sim}(k,r|p,t).\]
 The following proposition implies that $\mathbb{C}_{\sim}(k,r|p,t)\subseteq\mathbb{C}_{<}(k,r|p,t)$.
 \begin{prop}\label{defi-sim-1-5}
 For $1\leq j\leq 5$,  we have $\mathbb{C}^{(j)}_\sim(k,r|p,t)\subseteq\mathbb{C}_<(k,r|p,t)$. More precisely,
 \[\mathbb{C}^{(1)}_\sim(k,r|p,t)\subseteq\mathbb{C}^{(1)}_<(k,r|p,t),\]
 \[\mathbb{C}^{(2)}_\sim(k,r|p,t)\subseteq\mathbb{C}^{(1)}_<(k,r|p,t)\bigcup\mathbb{C}^{(2)}_<(k,r|p,t),\]
 \[\mathbb{C}^{(3)}_\sim(k,r|p,t)\subseteq\mathbb{C}^{(1)}_<(k,r|p,t)\bigcup\mathbb{C}^{(3)}_<(k,r|p,t),\]
 \[\mathbb{C}^{(4)}_\sim(k,r|p,t)\subseteq\mathbb{C}^{(1)}_<(k,r|p,t)\bigcup\mathbb{C}^{(2)}_<(k,r|p,t),\]
 and
\[\mathbb{C}^{(5)}_\sim(k,r|p,t)\subseteq\mathbb{C}^{(4)}_<(k,r|p,t)\bigcup\mathbb{C}^{(5)}_<(k,r|p,t).\]
\end{prop}

We will build a bijection $\mathcal{H}_{p,t}$ between $\mathbb{C}_{<}(k,r|p,t)$ and $\mathbb{C}_{\sim}(k,r|p,t)$ and build a bijection $\mathcal{I}_{p,t}$ between $\mathbb{C}_{\sim}(k,r|p,t)$ and $\mathbb{C}_{=}(k,r|p,t)$. Then, $\Phi_{p,t}=\mathcal{I}_{p,t}\cdot\mathcal{H}_{p,t}$ is a bijection between $\mathbb{C}_{<}(k,r|p,t)$ and $\mathbb{C}_{=}(k,r|p,t)$. Thus, Theorem \ref{equiv-main-lem-1-1} is equivalent to the following statements.

\begin{thm}\label{equiv-main-H}
There is a bijection $\mathcal{H}_{p,t}$ between $\mathbb{C}_{<}(k,r|p,t)$ and $\mathbb{C}_{\sim}(k,r|p,t)$. Moreover, for a partition $\pi\in \mathbb{C}_{<}(k,r|p,t)$, we have $\mu=\mathcal{H}_{p,t}(\pi)\in\mathbb{C}_{\sim}(k,r|p,t)$ such that
\begin{equation}\label{equiv-main-H-0-eqn}
|\mu|=|\pi|+2l \text{ and } \ell(\mu)=\ell(\pi),
\end{equation}
where $l$ is the largest integer such that $\pi^{(2)}_l>I_{p,t}(\pi)$.
\end{thm}

\begin{thm}\label{equiv-main-I}
There is a bijection $\mathcal{I}_{p,t}$ between $\mathbb{C}_{\sim}(k,r|p,t)$ and $\mathbb{C}_{=}(k,r|p,t)$. Moreover, for a partition $\mu\in \mathbb{C}_{\sim}(k,r|p,t)$, we have $\omega=\mathcal{I}_{p,t}(\mu)\in\mathbb{C}_{=}(k,r|p,t)$ such that
\begin{equation}\label{equiv-main-I-0-eqn}
|\omega|=|\mu|+2(p-l)+2t+1 \text{ and } \ell(\omega)=\ell(\mu)+1,
\end{equation}
where $l$ is the largest integer such that $\mu^{(2)}_l>I_{p,t}(\mu)$.
\end{thm}

\section{The dilation $\mathcal{H}_{p,t}$ and the reduction $\mathcal{R}_{p,t}$}

To give a proof of Theorem \ref{equiv-main-H}, we will introduce the dilation    $\mathcal{H}_{p,t}$ and its inverse the reduction $\mathcal{R}_{p,t}$, which are the modifications of the  dilation operation and the reduction operation respectively, introduced by He and Zhao \cite{He-Zhao-2023}. For more details, please see \cite[Section 5]{He-Zhao-2023}.

\subsection{The dilation    $\mathcal{H}_{p,t}$}

Let $\pi$ be a partition in $\mathbb{C}_{<}(k,r|p,t)$. Assume that $l$ is the largest integer such that $\pi^{(2)}_l>I_{p,t}(\pi)$. If $l\geq 1$, then the insertion types of $\pi^{(2)}_1,\pi^{(2)}_2,\ldots,$ $\pi^{(2)}_l$, introduced by He and Zhao \cite{He-Zhao-2023}, are defined as follows.

Set $b=0$ and $p_0=l+1$ and carry out the following procedure:

\begin{itemize}

    \item[(A)]  Assume that $p_{b+1}$ is the smallest integer such that $p_{b+1}\leq p_b-1$, $\pi^{(2)}_{p_{b+1}}-\pi^{(2)}_{p_b-1}=4(p_b-1-p_{b+1})$, and satisfying one of the following conditions:
    \begin{itemize}

\item[(1)] $\pi^{(2)}_{i}$ is of starting type $s_3$ and there exist parts $\pi^{(2)}_{i}$ with marks greater than $2$ in  $GG(\pi)$ for $p_{b+1}\leq i\leq p_b-1$. We say that $\pi^{(2)}_{p_b-1},\ldots,\pi^{(2)}_{p_{b+1}}$ are of insertion type $\check{A}_1${\rm;}

\item[(2)]  $\pi^{(2)}_{i}$ is of starting type $s_1$ for $p_{b+1}\leq i\leq p_b-1$. We say that  $\pi^{(2)}_{p_b-1},\ldots,\pi^{(2)}_{p_{b+1}}$ are of insertion type $\check{A}_2${\rm;}

\item[(3)] $\pi^{(2)}_{i}$ is of starting type $s_2$ for $p_{b+1}\leq i\leq p_b-1$ and  $\pi^{(2)}_{p_b-1}$ appears exactly once in $\pi$.
We say that  $\pi^{(2)}_{p_b-1},\ldots,\pi^{(2)}_{p_{b+1}}$ are of insertion type $\check{A}_3${\rm;}

\item[(4)] $\pi^{(2)}_{i}$ is of starting type $s_2$ and $\pi^{(2)}_{i}$  appears at least twice in $\pi$ for $p_{b+1}\leq i\leq p_b-1$.
We say that $\pi^{(2)}_{p_b-1},\ldots,\pi^{(2)}_{p_{b+1}}$ are of insertion type $\check{B}${\rm;}

\item[(5)] $\pi^{(2)}_{i}$ is of starting type $s_3$ for $p_{b+1}\leq i\leq p_b-1$ and   $\pi^{(2)}_{p_b-1}$ appears exactly twice in $\pi$. We say that  $\pi^{(2)}_{p_b-1},\ldots,\pi^{(2)}_{p_{b+1}}$ are of insertion type $\check{C}$.

\end{itemize}
Define $\check{\beta}_{p_b-1}(\pi)=\cdots=\check{\beta}_{p_{b+1}}(\pi)=\{p_{b+1},\ldots,p_b-1\}$.

\item[(B)] Replace $b$ by $b+1$. If $p_b=1$, then we are done. If $p_b>1$, then go back to (A).

\end{itemize}

For example, let $\pi$ be the partition in $\mathbb{C}_{<}(4,3|6,5)$ defined in \eqref{GG-example-1}. In such case, we have $I_{6,5}(\pi)=\pi^{(2)}_{5}=18$, and so $l=4$ is the largest integer such that $\pi^{(2)}_l=22>I_{6,5}(\pi)=18$. Then, we have
 \begin{itemize}
\item[(1)] $\pi^{(2)}_4=22$ is of insertion type $\check{A}_1$ and $\check{\beta}_{4}(\pi)=\{4\}$;

\item[(2)] $\pi^{(2)}_3=26$ is of insertion type $\check{C}$ and $\check{\beta}_{3}(\pi)=\{3\}$;

\item[(3)] $\pi^{(2)}_2=32$ and $\pi^{(2)}_1=36$ are of insertion type $\check{A}_3$ and $\check{\beta}_{2}(\pi)=\check{\beta}_{1}(\pi)=\{1,2\}$.
\end{itemize}

\begin{prop}\label{prop-ileq4-0}
Let $\pi$ be a partition in $\mathbb{C}_{<}(k,r|p,t)$ with $\pi^{(2)}_l=I_{p,t}(\pi)+2$ or $I_{p,t}(\pi)+4$, where $l$ is the largest integer such that $\pi^{(2)}_l>I_{p,t}(\pi)$.
\begin{itemize}
\item[\rm{(1)}] If $\pi\in\mathbb{C}^{(1)}_{<}(k,r|p,t)$, then  $\pi^{(2)}_l=2t+6$ is of insertion type $\check{A}_3$ or $\check{C}$.

\item[\rm{(2)}] If $\pi\in\mathbb{C}^{(2)}_{<}(k,r|p,t)$, then  $\pi^{(2)}_l=2t+6$ is of insertion type $\check{A}_1$ or $\check{A}_2$ or $\check{B}$.

\item[\rm{(3)}] If $\pi\in\mathbb{C}^{(3)}_{<}(k,r|p,t)$, then  $\pi^{(2)}_l=2t+6$ is of insertion type $\check{A}_1$.

\item[\rm{(4)}] If $\pi\in\mathbb{C}^{(4)}_{<}(k,r|p,t)$, then  $\pi^{(2)}_l=2t+4$ is of insertion type $\check{C}$.

\item[\rm{(5)}] If $\pi\in\mathbb{C}^{(5)}_{<}(k,r|p,t)$, then  $\pi^{(2)}_l=2t+4$ is of insertion type $\check{A}_1$.

\item[\rm{(6)}] For $6\leq j\leq 12$,   $\pi^{(2)}_l=I_{p,t}(\pi)+4$ is of insertion type $\check{A}_1$ or $\check{C}$ if $\pi\in\mathbb{C}^{(j)}_{<}(k,r|p,t)$.

\end{itemize}

\end{prop}

\pf By definition, we obtain that the conditions (1)-(5) hold. We proceed to show the condition (6). Assume that $\pi$ is a partition in $\mathbb{C}^{(j)}_{<}(k,r|p,t)$, where $6\leq j\leq 12$. Using the condition (3) in Proposition \ref{prop-i-1}, we see that $I_{p,t}(\pi)+2$ does not occur in $\pi$, and so $\pi^{(2)}_l=I_{p,t}(\pi)+4$. It follows from the definition of G\"ollnitz-Gordon marking that there is a $1$-marked $I_{p,t}(\pi)+4$ in $GG(\pi)$, which implies that $\pi^{(2)}_l$ is of insertion type $\check{A}_1$ or $\check{C}$. This completes the proof.  \qed

We need to recall definitions of special partition and the G\"ollnitz-Gordon marking of a special partition, introduced by He and Zhao \cite{He-Zhao-2023}. A special
partition $\pi$ is an ordinary partition in which the largest odd part in $\pi$ may be overlined.
The G\"ollnitz-Gordon marking of a special partition is given as follows.
\begin{defi}\label{defi-mark-special}
The G\"ollnitz-Gordon  marking  of a special partition $\pi$, denoted $\overline{GG}(\pi)$, is an assignment of positive integers {\rm(}marks{\rm)} to the parts of  $\pi=(\pi_1,\pi_2,\ldots,\pi_\ell)$ from  smallest to  largest such that the marks are as small as possible subject to the conditions that for $1\leq i\leq \ell$\rm{,}

\begin{itemize}
\item[{\rm (1)}] the integer assigned to $\pi_i$ is different from the integers assigned to the parts $\pi_g$ such that $\pi_i-\pi_g\leq 2$ with strict inequality if $\pi_i$ is an odd part for $g>i$\rm{;}

 \item[{\rm (2)}]  $\pi_i$ can not be assigned with $1$ if $\pi_i$ is an overlined odd part.

 \end{itemize}
\end{defi}

Now, we are in a position to give the definition of the dilation operation.

\begin{defi}
Let $\pi$ be a partition in $\mathbb{C}_{<}(k,r|p,t)$. Assume that $l$ is the largest integer such that $\pi^{(2)}_l>I_{p,t}(\pi)$. If $l=0$, then the dilation    $\mathcal{H}_{p,t}$ is defined as the identity map, that is,  $\mathcal{H}_{p,t}(\pi)=\pi$. If $l\geq 1$, then we define the dilation   $\mathcal{H}_{p,t}(\pi)$  as follows. There are three steps.
\begin{itemize}
	\item[] Step 1: We first do the operation related to $\pi^{(2)}_{l}$ and denote the resulting special partition by $\pi^l$, which is called the basic dilation of $\pi^{(2)}_{l}$. There are five cases.
	
		\begin{itemize}
	
		\item[] Case 1: $\pi^{(2)}_{l}$ is of insertion type $\check{A}_1$. We may write  $\pi^{(2)}_{l}=2t_{l}$. Let $r_{l}$ be the largest mark of parts  $2t_{l}$ in $GG(\pi)$. Then replace the $r_{l}$-marked $2t_{l}$ in $GG(\pi)$ by $2t_{l}+1$.
		
		\item[] Case 2: $\pi^{(2)}_{l}$ is of insertion type $\check{A}_2$. We may write  $\pi^{(2)}_{l}=2t_{l}+2$. Then set $r_{l}=1$ and replace the $1$-marked $2t_{l}$ in $GG(\pi)$ by $2t_{l}+1$.
		
			\item[] Case 3: $\pi^{(2)}_{l}$ is of insertion type $\check{A}_3$. We may write  $\pi^{(2)}_{l}=2t_{l}$. Then set $r_{l}=2$ and replace the $2$-marked $2t_{l}$ in $GG(\pi)$ by $\overline{2t_{l}+1}$.
			
				\item[] Case 4: $\pi^{(2)}_{l}$ is of insertion type $\check{B}$. We may write $\pi^{(2)}_{l}=2t_{l}$. Let $r_{l}$ be the largest mark of parts  $2t_{l}$ in ${GG}(\pi)$. Then replace the $r_{l}$-marked $2t_{l}$ in $GG(\pi)$ by $\overline{2t_{l}+1}$.
				
				\item[] Case 5: $\pi^{(2)}_{l}$ is of insertion type $\check{C}$. We may write $\pi^{(2)}_{l}=2t_{l}$. Then set $r_{l}=2$ and replace the $2$-marked $2t_{l}$ in $GG(\pi)$ by ${2t_{l}+1}$.
				
				\end{itemize}
				
		\item[] Step 2: If $l=1$, then go to Step 3 directly. If $l>1$, then set $b=l$ and repeat the following process.

		\begin{itemize}
			\item[] {\rm(A)}  There are the following two cases.
			
			\begin{itemize}
				\item[]  {\rm Case (A)-1:} $\check{\beta}_{b}(\pi)=\check{\beta}_{b-1}(\pi)$. In this case, we find that $2t_b+4$ occurs in $\pi^{b}$. Set  $t_{b-1}=t_b+2$ and $r_{b-1}$ to be the largest integer such that $r_{b-1}\leq r_b$ and there is an  $r_{b-1}$-marked $2t_b+4$ in $\overline{GG}(\pi^{b})$. Then $\pi^{b-1}$ is obtained by replacing the $r_b$-marked $2t_b+1$ {\rm(}resp. $\overline{2t_b+1}${\rm)} by $2t_b+2$ and replacing   the $r_{b-1}$-marked $2t_{b-1}$ in $\overline{GG}(\pi^{b})$ by  $2t_{b-1}+1$ {\rm(}resp. $\overline{2t_{b-1}+1}${\rm)}.

				\item[] {\rm Case (A)-2:} $\check{\beta}_{b}(\pi)\neq \check{\beta}_{b-1}(\pi)$. Then $\pi^{b-1}$ is obtained by replacing the $r_b$-marked $2t_b+1$ {\rm(}resp. $\overline{2t_b+1}${\rm)} in $\overline{GG}(\pi^{b})$ by $2t_b+2$ and doing the basic dilation of $\pi^{(2)}_{b-1}$.
				
			\end{itemize}
			
			\item[] {\rm(B)}  Replace $b$ by $b-1$. If $b>1$, then go back to {\rm(A)}. Otherwise, go to Step 3.
		\end{itemize}

	\item[] Step 3: There is an $r_1$-marked $2t_1+1$ {\rm(}resp. $\overline{2t_1+1}${\rm)} in $\overline{GG}(\pi^{1})$.  Then replace the  $r_1$-marked $2t_1+1$ {\rm(}resp. $\overline{2t_1+1}${\rm)} in $\overline{GG}(\pi^{1})$ by $2t_1+2$ and denote the resulting partition by $\mathcal{H}_{p,t}(\pi)$.
	
\end{itemize}
\end{defi}

For example, let $\pi$ be the partition in $\mathbb{C}_{<}(4,3|6,5)$ defined in \eqref{GG-example-1}. We apply dilation    $\mathcal{H}_{6,5}$  to $\pi$ to get a partition $\mu$. Here we also give the intermediate special partitions $\pi^4$, $\pi^3$, $\pi^2$ and $\pi^1$. The parts in boldface are the changed parts.

\[
\overline{GG}(\pi^4)=\left[
\begin{array}{cccccccccccccccccccc}
&&6&&&12&&16&&&{\bf 23}&&&&&&38\\
&2&6&&10&&14&&18&22&&26&&32&&36&\\
1&&6&9&&12&&16&&22&&26&30&&34&&38
\end{array}
\right].
\]
 \[\downarrow\]
 \[
\overline{GG}(\pi^3)=\left[
\begin{array}{cccccccccccccccccccc}
&&6&&&12&&16&&&{\bf 24}&&&&&&&38\\
&2&6&&10&&14&&18&22&&&{\bf 27}&&32&&36&\\
1&&6&9&&12&&16&&22&&26&&30&&34&&38
\end{array}
\right].
\]
\[\downarrow\]
 \[
\overline{GG}(\pi^2)=\left[
\begin{array}{cccccccccccccccccccc}
&&6&&&12&&16&&&24&&&&&&&38\\
&2&6&&10&&14&&18&22&&&{\bf 28}&&{\bf \overline{33}}&&36&\\
1&&6&9&&12&&16&&22&&26&&30&&34&&38
\end{array}
\right].
\]
\[\downarrow\]
 \[
\overline{GG}(\pi^1)=\left[
\begin{array}{cccccccccccccccccccc}
&&6&&&12&&16&&&24&&&&&&&38\\
&2&6&&10&&14&&18&22&&&28&&{\bf 34}&&{\bf \overline{37}}&\\
1&&6&9&&12&&16&&22&&26&&30&34&&&38
\end{array}
\right].
\]
\[\downarrow\]
 \begin{equation}\label{example-H-R}
{GG}(\mu)=\left[
\begin{array}{cccccccccccccccccccc}
&&6&&&12&&16&&&24&&&&&38\\
&2&6&&10&&14&&18&22&&&28&&34&{\bf 38}\\
1&&6&9&&12&&16&&22&&26&&30&34&38
\end{array}
\right].
\end{equation}

With a similar argument as in  \cite[Section 5.4]{He-Zhao-2023}, we get the following lemma, which says that the dilation $\mathcal{H}_{p,t}$ is a map from $\mathbb{C}_{<}(k,r|p,t)$ to $\mathbb{C}_{\sim}(k,r|p,t)$.
\begin{lem}\label{equiv-main-H-lem-1}
For $1\leq j\leq 12$, the dilation $\mathcal{H}_{p,t}$ is a map from $\mathbb{C}^{(j)}_{<}(k,r|p,t)$ to $\mathbb{C}^{(j)}_{\sim}(k,r|p,t)$. Moreover, for a partition $\pi\in \mathbb{C}^{(j)}_{<}(k,r|p,t)$, we have $\mu=\mathcal{H}_{p,t}(\pi)\in\mathbb{C}^{(j)}_{\sim}(k,r|p,t)$ such that
\[|\mu|=|\pi|+2l \text{ and } \ell(\mu)=\ell(\pi),\]
where $l$ is the largest integer such that $\pi^{(2)}_l>I_{p,t}(\pi)$.
\end{lem}

\subsection{The reduction $\mathcal{R}_{p,t}$}

In this subsection, we introduce the reduction $\mathcal{R}_{p,t}$, which will be shown to be inverse map of the dilation $\mathcal{H}_{p,t}$.

\begin{defi}
Let $\mu$ be a partition in $\mathbb{C}_{\sim}(k,r|p,t)$. Assume that $l$ is the largest integer such that $\mu^{(2)}_l>I_{p,t}(\mu)$. If $l=0$, then the reduction    $\mathcal{R}_{p,t}$ is defined as the identity map, that is,  $\mathcal{R}_{p,t}(\mu)=\mu$. If $l\geq 1$, then we define the reduction  $\mathcal{R}_{p,t}(\pi)$  as follows. There are three steps.

\begin{itemize}

\item[] Step 1: We do the following operation related to $\mu^{(2)}_1$, called the basic reduction of $\mu^{(2)}_1$. There are  five cases.

          \begin{itemize}
              \item[] Case 1: $\mu^{(2)}_1$ is of reduction type  $\hat{A}_1$.   We may write         $\mu^{(2)}_1=2t_1$.  Then, there is a $1$-marked part $2t_1$ in ${GG}(\mu)$ and there exist parts $2t_1+2$ in $\mu$. Let $r_1$ be the smallest mark of parts  $2t_1+2$ in ${GG}(\mu)$.  Then replace the $r_1$-marked $2t_1+2$ in ${GG}(\mu)$ by $2t_1+1$ to get $\mu^1$.

                  \item[] Case 2: $\mu^{(2)}_1$ is of reduction type  $\hat{A}_2$. We may write $\mu^{(2)}_1=2t_1+2$. Then, there is a $1$-marked part $2t_1+2$ in ${GG}(\mu)$ and there do not exist parts $2t_1+4$ in $\mu$. Set $r_1=1$ and replace the $1$-marked $2t_1+2$ in ${GG}(\mu)$ by $2t_1+1$ to obtain $\mu^1$.

                        \item[] Case 3: $\mu^{(2)}_1$ is of reduction type  $\hat{A}_3$. We may write $\mu^{(2)}_1=2t_1+2$. Then, there are $1$-marked parts $2t_1-2$ and $2t_1+2$ in ${GG}(\mu)$ and there do not exist parts $2t_1+4$ in $\mu$. Set $r_1=2$ and replace the $2$-marked $2t_1+2$ in ${GG}(\mu)$ by $\overline{2t_1+1}$ to obtain $\mu^1$.

                      \item[] Case 4:  $\mu^{(2)}_1$ is of reduction type  $\hat{B}$. We may write $\mu^{(2)}_1=2t_1$.  Then,  there exist parts $2t_1+2$ with mark $1$ and  marks  greater than $2$ in ${GG}(\mu)$. Let $r_1$ be the smallest mark except for $1$ of parts  $2t_1+2$ in ${GG}(\mu)$.  Then  $\mu^1$ is obtained by replacing the $r_1$-marked $2t_1+2$ in ${GG}(\mu)$ by $\overline{2t_1+1}$.

                     \item[] Case 5: $\mu^{(2)}_1$ is of reduction type  $\hat{C}$. We may write $\mu^{(2)}_1=2t_1+2$. Then, there is a $1$-marked part $2t_1$   and there do not exist parts $2t_1+4$ with marks greater than $2$ in ${GG}(\mu)$. Set $r_1=2$ and replace the $2$-marked $2t_1+2$ in ${GG}(\mu)$ by $2t_1+1$ to get $\mu^1$.

          \end{itemize}

          \item[] Step 2: If $l=1$,  then go to Step 3 directly. If $l>1$, then set $b=1$ and repeat the following process.

\begin{itemize}
    \item[{\rm(A)}]  Replace the $r_{b}$-marked $2t_{b}+1$ {\rm(}resp. $\overline{2t_{b}+1}${\rm)} in $\overline{GG}(\mu^{b})$ by $2t_{b}$ and apply the basic reduction of $\mu^{(2)}_{b+1}$    to get $\mu^{b+1}$.

     \item[{\rm(B)}] replace $b$ by $b+1$. If $b<l$, then go back to {\rm(A)}. Otherwise, go to Step 3.
\end{itemize}

\item[] Step 3: There is an $r_l$-marked $2t_l+1$ {\rm(}resp. $\overline{2t_l+1}${\rm)} in $\overline{GG}(\mu^{l})$.  Then replace the  $r_l$-marked $2t_l+1$ {\rm(}resp. $\overline{2t_l+1}${\rm)} in $\overline{GG}(\mu^{l})$ by $2t_l$ and denote the resulting partition by $\mathcal{R}_{p,t}(\mu)$.

\end{itemize}

\end{defi}

For example, let $\mu$ be the partition in $\mathbb{C}_{\sim}(4,3|6,5)$, whose G\"ollnitz-Gordon marking is given in \eqref{example-H-R}. Applying the reduction $\mathcal{R}_{6,5}$, then the same process to get $\mu$ could be run in reverse. Then, we can obtain the partition $\pi=\mathcal{R}_{6,5}(\mu)$, which is the partition with G\"ollnitz-Gordon marking give in \eqref{GG-example-1}.

With a similar argument as in \cite[Section 5.4]{He-Zhao-2023}, we get the following lemma, which says that the reduction $\mathcal{R}_{p,t}$ is a map from $\mathbb{C}_{\sim}(k,r|p,t)$ to $\mathbb{C}_{<}(k,r|p,t)$.
\begin{lem}\label{equiv-main-R-lem-1}
For $1\leq j\leq 12$, the reduction $\mathcal{R}_{p,t}$ is a map from $\mathbb{C}^{(j)}_{\sim}(k,r|p,t)$ to $\mathbb{C}^{(j)}_{<}(k,r|p,t)$. Moreover, for a partition $\mu\in \mathbb{C}^{(j)}_{\sim}(k,r|p,t)$, we have $\pi=\mathcal{R}_{p,t}(\mu)\in\mathbb{C}^{(j)}_{<}(k,r|p,t)$ such that
\[|\pi|=|\mu|-2l \text{ and } \ell(\pi)=\ell(\mu),\]
where $l$ is the largest integer such that $\mu^{(2)}_l>I_{p,t}(\mu)$.

\end{lem}

\subsection{Proof of Theorem \ref{equiv-main-H}}

We are now in a position to give a proof of Theorem \ref{equiv-main-H}.

{\noindent \bf Proof of Theorem \ref{equiv-main-H}.} Using Lemma \ref{equiv-main-H-lem-1}, we know that the dilation $\mathcal{H}_{p,t}$ is a map from $\mathbb{C}_{<}(k,r|p,t)$ to $\mathbb{C}_{\sim}(k,r|p,t)$ satisfying \eqref{equiv-main-H-0-eqn}. Appealing to Lemma \ref{equiv-main-R-lem-1}, we deduce that the reduction $\mathcal{R}_{p,t}$ is a map from $\mathbb{C}_{\sim}(k,r|p,t)$ to $\mathbb{C}_{<}(k,r|p,t)$. With a similar argument as in  \cite[Section 5.4]{He-Zhao-2023}, we find that the dilation $\mathcal{H}_{p,t}$ and the reduction $\mathcal{R}_{p,t}$ are inverses of each other. The proof is complete.   \qed

\section{The insertion  $\mathcal{I}_{p,t}$ and the separation  $\mathcal{S}_{p,t}$}

In this section, we will introduce the insertion  $\mathcal{I}_{p,t}$ and the separation  $\mathcal{S}_{p,t}$ and then give a proof of Theorem \ref{equiv-main-I}.

\subsection{The insertion  $\mathcal{I}_{p,t}$}
 We will define the insertion   $\mathcal{I}_{p,t}$ from $\mathbb{C}_{\sim}(k,r|p,t)$ to $\mathbb{C}_{=}(k,r|p,t)$ in this subsection.  Before doing this, we give the following  lemma, which plays an important role in considering the mark of $2t+1$ in the resulting partition.

 \begin{prop}\label{prop-i-3}

Let $\mu$ be a partition in $\mathbb{C}_{\sim}(k,r|p,t)$ such that there is a $1$-marked $2t$ in $GG(\mu)$. Then, there is no $2$-marked $2t$ in $GG(\mu)$.
 	
 \end{prop}

 \pf Suppose to the contrary that there is a $2$-marked $2t$ in $GG(\mu)$. Note that $\mathbb{C}_{\sim}(k,r|p,t)\subseteq\mathbb{C}_{<}(k,r|p,t)$, then by the condition (2) in Definition \ref{defi-c<}, we see that $\mu^{(2)}_{p+1}<2t+1<\mu^{(2)}_{p}$. It yields $\mu^{(2)}_{p+1}=2t$. Under the condition that  is a $1$-marked $2t$ in $GG(\mu)$, we see that $\mu^{(2)}_{p+1}=2t$ is of starting type $s_3$, which contradicts the condition (4) in Definition \ref{defi-c<}. This completes the proof.  \qed

 To give the insertion   $\mathcal{I}_{p,t}$,  we will define the $j$-th kind of the insertion   $\mathcal{I}^{(j)}_{p,t}$ from $\mathbb{C}^{(j)}_{\sim}(k,r|p,t)$ to $\mathbb{C}^{(j)}_{=}(k,r|p,t)$ for $1\leq j\leq 12$. We first state the  $j$-th kind of the insertion   $\mathcal{I}^{(j)}_{p,t}$ for $1\leq j\leq 5$.

 \begin{defi}\label{defi-insertion-1-5}
 	
 For $1\leq j\leq 5$,	let $\mu$ be a partition in $\mathbb{C}^{(j)}_{\sim}(k,r|p,t)$.  Define the $j$-th kind of the insertion   $\mathcal{I}^{(j)}_{p,t}$ as follows{\rm:} add $2t+1$ as a part of $\mu$.

 \end{defi}

 The following lemma says that $\mathcal{I}^{(j)}_{p,t}$ is a map from $\mathbb{C}^{(j)}_{\sim}(k,r|p,t)$ to $\mathbb{C}^{(j)}_{=}(k,r|p,t)$ for $1\leq j\leq 5$.
\begin{lem}\label{lem-insertion-1-5}
For $1\leq j\leq 5$, let $\mu$ be a partition in $\mathbb{C}^{(j)}_{\sim}(k,r|p,t)$ and let $\omega=\mathcal{I}^{(j)}_{p,t}(\mu)$. Then, $\omega$ is a partition in $\mathbb{C}^{(j)}_{=}(k,r|p,t)$ such that
\[|\omega|=|\mu|+2t+1\text{ and }\ell(\omega)=\ell(\mu)+1.\]
\end{lem}

\pf By the construction of $\omega$, we find that the largest odd part of $\omega$ is $2t+1$, $|\omega|=|\mu|+2t+1$, $\ell(\omega)=\ell(\mu)+1$, and the marks of parts not exceeding $2t$ in $GG(\omega)$ are the same as those in $GG(\mu)$. Using Proposition \ref{prop-i-3}, we see that the mark of $2t+1$ in $GG(\omega)$ is at most $2$. More precisely, if there is no $1$-marked $2t$ in $GG(\mu)$, then $2t+1$ is marked with $1$ in  $GG(\omega)$; if there is a $1$-marked $2t$ in $GG(\mu)$, then $2t+1$ is marked with $2$ in  $GG(\omega)$.

It follows from Propositions \ref{prop-i-1} and \ref{defi-sim-1-5} that
$2t+2$ and $2t+4$ can not both occur in $\mu$, and so the marks of parts greater than or equal to $2t+4$ in $GG(\omega)$ are the same as those in $GG(\mu)$. It yields that $\omega^{(2)}_p=\mu^{(2)}_p\geq 2t+4$ and $\omega^{(2)}_{p+1}\leq 2t+2$.

Assume that $2t+2$ occurs in $\omega$, then by the construction of $\omega$, we deduce that $2t+2$ occurs in $\mu$, and so $\mu\in\mathbb{C}^{(1)}_{\sim}(k,r|p,t)$ or $\mathbb{C}^{(3)}_{\sim}(k,r|p,t)$. In view of Propositions \ref{prop-i-1} and \ref{defi-sim-1-5}, we know that $2t+2$ occurs once in $\mu$, $2t+2$ is marked with $1$ in $GG(\mu)$ and $2t+4$ does not occur in $\mu$, and so there is no $1$-marked $2t$ in $GG(\mu)$ and $2t+4$ does not occur in $\omega$.
Moreover, there is no $2$-marked $2t$ in $GG(\mu)$, otherwise we obtain that $\mu^{(2)}_{p+1}=2t$ and it is of starting type $s_2$, which contradicts the condition (4) in Definition \ref{defi-c<}. Therefore, $2t+1$ is marked with $1$ in $GG(\omega)$ and $2t+2$ is marked with $2$ in $GG(\omega)$. Note that $2t+4$ does not occur in $\omega$, we obtain that $\omega^{(2)}_{p+1}=2t+2$ and it is of starting type $s_0$.

Now, we conclude that  $\omega$ is a partition in $\mathbb{C}_{=}(k,r|p,t)$. Then, we proceed to show that for $1\leq j\leq 5$, $\omega$ is a partition in $\mathbb{C}^{(j)}_{=}(k,r|p,t)$. We consider the following five cases.

Case 1: $j=1$. In this case,  we have  $\omega^{(2)}_p=\mu^{(2)}_p\geq 2t+8$, and so $\omega\in\mathbb{C}^{(1)}_{=}(k,r|p,t)$.

Case 2: $j=2$. Since $\mu\in\mathbb{C}^{(2)}_{\sim}(k,r|p,t)$, we know that $2t+2$ does not occur in $\mu$ and $\mu^{(2)}_{p}=2t+6$ is of reduction type $\hat{A}_1$ or $\hat{A}_2$ or $\hat{B}$.
By the construction of $\omega$, we deduce that  $2t+2$ also does not occur in $\omega$, which implies that $\omega^{(2)}_{p+1}<2t+2$.

Assume that $\mu^{(2)}_{p}=2t+6$ is of reduction type $\hat{A}_1$ or $\hat{A}_2$, then $\omega^{(2)}_{p}=2t+6$ is of starting type $s_3$. Assume that $\mu^{(2)}_{p}=2t+6$ is of reduction type $\hat{B}$, then $\omega^{(2)}_{p}=2t+6$ is of starting type $s_2$, and $\omega^{(2)}_{i}+2$ occurs at least twice in $\omega$ if $\omega^{(2)}_{i}=\omega^{(2)}_{p}+4(p-i)$ for  $i\leq p$. So, we have $\omega\in\mathbb{C}^{(2)}_{=}(k,r|p,t)$.

Case 3: $j=3$. It follows from $\mu\in\mathbb{C}^{(3)}_{\sim}(k,r|p,t)$ that $2t+2$ occurs in $\mu$ and
 $\mu^{(2)}_{p}=2t+6$ is of reduction type $\hat{A}_1$. From the proof above, we have $\omega^{(2)}_{p+1}=2t+2$.
By the construction of $\omega$, we see that  $\omega^{(2)}_{p}=2t+6$ is of starting type $s_3$, and if
$\omega^{(2)}_i=\omega^{(2)}_{p}+4(p-i)$ then  $\omega^{(2)}_i+2$ occurs in $\omega$ for $i\leq p$. We arrive at $\omega\in\mathbb{C}^{(3)}_{=}(k,r|p,t)$.

Case 4: $j=4$. In this case, we see that $\mu^{(2)}_{p}=2t+6$ is of reduction type $\hat{C}$.  Then, $\omega^{(2)}_{p}=2t+6$ and it is of starting type $s_1$ or $s_2$, and if $\omega^{(2)}_{p}$ is of starting type $s_2$ then there exists $i$ such that $i\leq p$, $\omega^{(2)}_{i}=\omega^{(2)}_{p}+4(p-i)$ and  $\omega^{(2)}_{i}+2$ occurs once in $\omega$.
With a similar argument as in Case 2, we get $\omega^{(2)}_{p+1}<2t+2$, which yields $\omega\in\mathbb{C}^{(4)}_{=}(k,r|p,t)$.

Case 5: $j=5$. It is immediate from the construction of $\omega$ that $\omega\in\mathbb{C}^{(5)}_{=}(k,r|p,t)$. Thus, we have completed the proof.  \qed

Next, we give the $j$-th kind of the insertion $\mathcal{I}^{(j)}_{p,t}$ for $6\leq j\leq 12$.

\begin{defi}\label{defi-insertion-6}
		Let $\mu$ be a partition in $\mathbb{C}^{(6)}_{\sim}(k,r|p,t)$. Assume that $p_1$ is the first starting cluster index of $\mu$. Define $\mathcal{I}^{(6)}_{p,t}\colon \mu \rightarrow \omega$ as follows{\rm:}
	add  $2t+1$ as a $1$-marked part into $GG(\mu)$ and replace the $1$-marked parts $\mu^{(2)}_p-2,\ldots,\mu^{(2)}_{p_1}-2$ in $GG(\mu)$ by $1$-marked parts $\mu^{(2)}_p,\ldots,\mu^{(2)}_{p_1}$ respectively to get $\omega$.
\end{defi}

For example, let $\mu$ be a partition in $\mathbb{C}^{(6)}_{\sim}(4,3|2,1)$ with G\"ollnitz-Gordon marking
\begin{equation}\label{example-<-6}
{GG}(\mu)=\left[
\begin{array}{cccccccccccccccccccc}
&&4&&8\\
&2&&6&&10\\
1&&4&&8
\end{array}
\right].
\end{equation}
It can be checked that $p_1=1$. Adding $3$ as a $1$-marked part into $GG(\mu)$ and replacing the $1$-marked parts $4$ and $8$ in $GG(\mu)$ by $1$-marked parts $6$ and $10$, we get
\begin{equation}\label{example-=-6}
GG(\omega)=\left[
\begin{array}{cccccccccccccccccccc}
&&&4&&8\\
&2&&&6&&10\\
1&&3&&6&&10
\end{array}
\right].
\end{equation}

\begin{defi}\label{defi-insertion-7}
		Let $\mu$ be a partition in $\mathbb{C}^{(7)}_{\sim}(k,r|p,t)$. Assume that $p_1$ is the first starting cluster index of $\mu$. Define  $\mathcal{I}^{(7)}_{p,t}\colon \mu \rightarrow \omega$ as follows{\rm:}

	{\rm(1)} Add  $2t+1$ as a part into $\mu$ and denote the resulting partition by $\nu$. Moreover, $2t+1$ is marked with $1$ in $GG(\nu)$, the parts $\mu^{(2)}_p-2,\mu^{(2)}_p+2,\ldots,\mu^{(2)}_{p_1}+2$ marked with $1$ in  $GG(\mu)$ are marked with $2$ in $GG(\nu)$, and the parts $\mu^{(2)}_p,\ldots,\mu^{(2)}_{p_1}$ marked with $2$ in  $GG(\mu)$ are marked with $1$ in $GG(\nu)$.
	
	{\rm(2)} Let $r_p$ be the largest mark of parts $\mu^{(2)}_p$ in  $GG(\nu)$. For $p_1\leq i<p$, assume that $r_{i+1}$ has been defined, then $r_i$ is defined to be the largest integer such that $r_i\leq r_{i+1}$ and there is an $r_i$-marked $\mu^{(2)}_i$ in $GG(\nu)$. Replace the $r_i$-marked  $\mu^{(2)}_i$ in $GG(\nu)$ by  $r_i$-marked  $\mu^{(2)}_i+2$ for $p_1\leq i\leq p$ to get $\omega$.
	
	\end{defi}

For example, let $\mu$ be a partition in $\mathbb{C}^{(7)}_{\sim}(4,3|3,1)$ with G\"ollnitz-Gordon marking
\begin{equation}\label{example-<-7}
{GG}(\mu)=\left[
\begin{array}{cccccccccccccccccccc}
&&6&&&&14\\
&&6&&10&&14\\
1&4&&8&&12&&16
\end{array}
\right].
\end{equation}
It can be checked that $p_1=1$. We add $3$ as a part into $\mu$ to get
\begin{equation}\label{example-<=-7}
{GG}(\nu)=\left[
\begin{array}{cccccccccccccccccccc}
&&&6&&&&14\\
&&4&&8&&12&&16\\
1&3&&6&&10&&14&
\end{array}
\right].
\end{equation}
Moreover, we have $r_3=3$, $r_2=1$ and $r_1=1$. Then, $\omega$ is obtained by replacing the $3$-marked $6$, $1$-marked $10$ and $1$-marked $14$ in ${GG}(\nu)$ by $3$-marked $8$, $1$-marked $12$ and $1$-marked $16$ respectively.
\begin{equation}\label{example-=-7}
GG(\omega)=\left[
\begin{array}{cccccccccccccccccccc}
&&&&8&&&14\\
&&4&&8&&12&&16\\
1&3&&6&&&12&&16
\end{array}
\right].
\end{equation}

	\begin{defi}\label{defi-insertion-8}
			Let $\mu$ be a partition in $\mathbb{C}^{(8)}_{\sim}(k,r|p,t)$. Assume that $p_1$ is the first starting cluster index of $\mu$. Define $\mathcal{I}^{(8)}_{p,t}\colon \mu \rightarrow \omega$ as follows{\rm:}
		add  $2t+1$ as a $2$-marked part into $GG(\mu)$ and replace the $2$-marked parts $\mu^{(2)}_p,\ldots,\mu^{(2)}_{p_1}$ in $GG(\mu)$ by $2$-marked parts $\mu^{(2)}_p+2,\ldots,\mu^{(2)}_{p_1}+2$ respectively to get $\omega$.
		\end{defi}

For example, let $\mu$ be a partition in $\mathbb{C}^{(8)}_{\sim}(4,3|3,2)$, whose G\"ollnitz-Gordon marking reads
\begin{equation}\label{example-<-8}
{GG}(\mu)=\left[
\begin{array}{cccccccccccccccccccc}
&&&6&&10&&&16\\
&2&&6&&10&&14\\
1&&4&&8&&12&&16
\end{array}
\right].
\end{equation}
It can be checked that $p_1=1$. Add  $5$ as a $2$-marked part into $GG(\mu)$ and replace the $2$-marked $6$, $10$ and $14$  in $GG(\mu)$ by $2$-marked $8$, $12$ and $16$ respectively. We get
\begin{equation}\label{example-=-8}
{GG}(\omega)=\left[
\begin{array}{cccccccccccccccccccc}
&&&&6&&10&&&16\\
&2&&5&&8&&12&&16\\
1&&4&&&8&&12&&16
\end{array}
\right].
\end{equation}

\begin{defi}\label{defi-insertion-9}
Let $\mu$ be a partition in $\mathbb{C}^{(9)}_{\sim}(k,r|p,t)$. Assume that $p_1$ is the first starting cluster index of $\mu$. Define $\mathcal{I}^{(9)}_{p,t}\colon \mu \rightarrow \omega$ as follows{\rm:}
add  $2t+1$ as a $1$-marked part into $GG(\mu)$ and replace the $1$-marked parts $\mu^{(2)}_p,\ldots,\mu^{(2)}_{p_1}$ in $GG(\mu)$ by $1$-marked parts $\mu^{(2)}_p+2,\ldots,\mu^{(2)}_{p_1}+2$ respectively to get $\omega$.
\end{defi}

For example, let $\mu$ be a partition in $\mathbb{C}^{(9)}_{\sim}(4,3|3,0)$, whose G\"ollnitz-Gordon marking is
\begin{equation}\label{example-<-9}
{GG}(\mu)=\left[
\begin{array}{cccccccccccccccccccc}
&4&&&12\\
2&&6&10\\
2&&6&10
\end{array}
\right].
\end{equation}
It can be checked that $p_1=1$. Adding  $1$ as a $1$-marked part into $GG(\mu)$ and replacing the $1$-marked $2$, $6$ and $10$  in $GG(\mu)$ by $1$-marked $4$, $8$ and $12$ respectively, we get
\begin{equation}\label{example-=-9}
{GG}(\omega)=\left[
\begin{array}{cccccccccccccccccccc}
&&4&&&&12\\
&2&&6&&10\\
1&&4&&8&&12
\end{array}
\right].
\end{equation}

\begin{defi}\label{defi-insertion-10}
Let $\mu$ be a partition in $\mathbb{C}^{(10)}_{\sim}(k,r|p,t)$. Assume that $p_1$ is the first starting cluster index of $\mu$. Define $\mathcal{I}^{(10)}_{p,t}\colon \mu \rightarrow \omega$ as follows{\rm:}
add $2t+1$ as a $1$-marked part into $GG(\mu)$ and replace the $1$-marked parts $\mu^{(2)}_p,\ldots,\mu^{(2)}_{p_1}$ in $GG(\mu)$ by $1$-marked parts $\mu^{(2)}_p+2,\ldots,\mu^{(2)}_{p_1}+2$ respectively to get $\omega$. Then, the  part $\mu^{(2)}_{p_1}+4$ marked with $1$ in  $GG(\mu)$ is  marked with   $2$ in $GG(\omega)$.
\end{defi}

For example, let $\mu$ be a partition in $\mathbb{C}^{(10)}_{\sim}(4,3|3,0)$ with G\"ollnitz-Gordon marking
\begin{equation}\label{example-<-10}
{GG}(\mu)=\left[
\begin{array}{cccccccccccccccccccc}
&4&&&12&\\
2&&6&10&&\\
2&&6&10&&14
\end{array}
\right].
\end{equation}
We find that $p_1=1$. Add  $1$ as a $1$-marked part into $GG(\mu)$ and replace the $1$-marked $2$, $6$ and $10$ in $GG(\mu)$ by $1$-marked $4$, $8$ and $12$  respectively to get $\omega$. Then, the  part $14$ marked with $1$ in  $GG(\mu)$ is  marked with   $2$ in $GG(\omega)$. We have
\begin{equation}\label{example-=-10}
{GG}(\omega)=\left[
\begin{array}{cccccccccccccccccccc}
&&4&&&&12&\\
&2&&6&&10&&14\\
1&&4&&8&&12&
\end{array}
\right].
\end{equation}

\begin{defi}\label{defi-insertion-11}
Let $\mu$ be a partition in $\mathbb{C}^{(11)}_{\sim}(k,r|p,t)$.  Assume that $p_1$ and $p_2$ are the first and the second starting cluster indexes of $\mu$ respectively.
Define $\mathcal{I}^{(11)}_{p,t}\colon \mu \rightarrow \omega$ as follows{\rm:}
add  $2t+1$ as a $1$-marked part into $GG(\mu)$ and replace the $1$-marked parts $\mu^{(2)}_p,\ldots,\mu^{(2)}_{p_1},\mu^{(2)}_{p_1-1}-2,\ldots,\mu^{(2)}_{p_2}-2$ in $GG(\mu)$ by $1$-marked parts $\mu^{(2)}_p+2,\ldots,\mu^{(2)}_{p_1}+2,\mu^{(2)}_{p_1-1},\ldots,\mu^{(2)}_{p_2}$ respectively to get $\omega$.
\end{defi}

For example, let $\mu$ be a partition in $\mathbb{C}^{(11)}_{\sim}(4,3|5,0)$, whose G\"ollnitz-Gordon marking is
\begin{equation}\label{example-<-11}
{GG}(\mu)=\left[
\begin{array}{cccccccccccccccccccc}
&4&&&12&&16&&\\
2&&6&10&&&16&&20\\
2&&6&10&&14&&18
\end{array}
\right].
\end{equation}
It can be checked that $p_1=3$ and $p_2=1$. Adding  $1$ as a $1$-marked part into $GG(\mu)$ and replacing the $1$-marked $2$, $6$, $10$, $14$ and $18$ in $GG(\mu)$ by $1$-marked $4$, $8$, $12$, $16$ and $20$  respectively, we get
\begin{equation}\label{example-=-11}
{GG}(\omega)=\left[
\begin{array}{cccccccccccccccccccc}
&&4&&&&12&&16&&\\
&2&&6&&10&&&16&&20\\
1&&4&&8&&12&&16&&20
\end{array}
\right].
\end{equation}

\begin{defi}\label{defi-insertion-12}
								
Let $\mu$ be a partition in $\mathbb{C}^{(12)}_{\sim}(k,r|p,t)$. Assume that $p_1$ and $p_2$ are the first and the second starting cluster indexes of $\mu$ respectively. Define $\mathcal{I}^{(12)}_{p,t}\colon \mu \rightarrow \omega$ as follows{\rm:}
								
{\rm(1)} Add $2t+1$ as a $1$-marked part into $GG(\mu)$ and replace the $1$-marked parts $\mu^{(2)}_p,\ldots,\mu^{(2)}_{p_1}$ in $GG(\mu)$ by $1$-marked parts $\mu^{(2)}_p+2,\ldots,\mu^{(2)}_{p_1}+2$ respectively. Denote the resulting partition by $\nu$. Then, the parts $\mu^{(2)}_{p_1-1}-2,\mu^{(2)}_{p_1-1}+2,\ldots,\mu^{(2)}_{p_2}+2$ marked with $1$ in  $GG(\mu)$ are marked with $2$ in $GG(\nu)$, and the parts $\mu^{(2)}_{p_1-1},\ldots,\mu^{(2)}_{p_2}$ marked with $2$ in  $GG(\mu)$ are marked with $1$ in $GG(\nu)$.

{\rm(2)} Let $r_{p_1-1}$ be the largest mark of parts $\mu^{(2)}_{p_1-1}$ in  $GG(\nu)$. For $p_2\leq i<p_1-1$, assume that $r_{i+1}$ has been defined, then $r_i$ is defined to be the largest integer such that $r_i\leq r_{i+1}$ and there is a $r_i$-marked $\mu^{(2)}_i$ in $GG(\nu)$. Repalce the $r_i$-marked  $\mu^{(2)}_i$ in $GG(\nu)$ by  $r_i$-marked  $\mu^{(2)}_i+2$ for $p_2\leq i\leq p_1-1$ to get $\omega$.
\end{defi}

For example, let $\mu$ be a partition in $\mathbb{C}^{(12)}_{\sim}(4,3|6,0)$ with G\"ollnitz-Gordon marking
\begin{equation}\label{example-<-12}
{GG}(\pi)=\left[
\begin{array}{cccccccccccccccccccc}
&4&&&12&&16&&&&24\\
2&&6&10&&&16&&20&&24&\\
2&&6&10&&14&&18&&22&&26
\end{array}
\right].
\end{equation}
It can be checked that $p_1=4$ and $p_2=1$. Add $1$ as a $1$-marked part into $GG(\mu)$ and replace the $1$-marked $2$, $6$ and $10$ in $GG(\mu)$ by $1$-marked $4$, $8$ and $12$  respectively to get
\begin{equation}\label{example-<=-12}
{GG}(\nu)=\left[
\begin{array}{cccccccccccccccccccc}
&&4&&&&12&&16&&&&24\\
&2&&6&&10&&14&&18&&22&&26\\
1&&4&&8&&12&&16&&20&&24&
\end{array}
\right].
\end{equation}
Moreover, we have $r_3=3$, $r_2=1$ and $r_1=1$. Then, $\omega$ is obtained by replacing the $3$-marked $16$, $1$-marked $20$ and $1$-marked $24$ in ${GG}(\nu)$ by $3$-marked $18$, $1$-marked $22$ and $1$-marked $26$ respectively.
\begin{equation}\label{example-=-12}
{GG}(\omega)=\left[
\begin{array}{cccccccccccccccccccc}
&&4&&&&12&&&18&&24\\
&2&&6&&10&&14&&18&22&&26\\
1&&4&&8&&12&&16&&22&&26
\end{array}
\right].
\end{equation}

Then, we proceed to show that for $6\leq j\leq 12$, the $j$-th kind of the insertion $\mathcal{I}^{(j)}_{p,t}$ is a map from $\mathbb{C}^{(j)}_{\sim}(k,r|p,t)$ to $\mathbb{C}^{(j)}_{=}(k,r|p,t)$.	

\begin{lem}\label{lem-insertion-6-12}
For $6\leq j\leq 12$, let $\mu$ be a partition in $\mathbb{C}^{(j)}_{\sim}(k,r|p,t)$ and let $\omega=\mathcal{I}^{(j)}_{p,t}(\mu)$. Then, $\omega$ is a partition in $\mathbb{C}^{(j)}_{=}(k,r|p,t)$ such that $D_{p,t}(\omega)=I_{p,t}(\mu),$
\[|\omega|=|\mu|+2(p-l)+2t+1\text{ and }\ell(\omega)=\ell(\mu)+1,\]
where $l$ is the largest integer such that $\mu^{(2)}_{l}>I_{p,t}(\mu)$.
\end{lem}

\pf For $6\leq j\leq 12$, it is clear from the definition of $\mathcal{I}^{(j)}_{p,t}$ that
\begin{itemize}			
\item[\rm{(1)}] the $j$-th kind of the insertion $\mathcal{I}^{(j)}_{p,t}$ is well-defined{\rm;}
			
\item[\rm{(2)}] the largest odd part of $\omega$ is $2t+1$ and the mark of $2t+1$ in $GG(\omega)$ is at  most $2${\rm;}
		
\item[\rm{(3)}] $\omega^{(2)}_p\geq \mu^{(2)}_p\geq 2t+2$ and $\omega^{(2)}_{p+1}\leq 2t+2${\rm;}

\item[\rm{(4)}] $I_{p,t}(\mu)+2$ does not occur in $\omega${\rm;}
		
\item[\rm{(5)}] $\omega^{(2)}_{i}$ is of starting type $s_3$ and $\omega^{(2)}_{i}+2$ occurs in $\omega$ if $\omega^{(2)}_{i}=I_{p,t}(\mu)+4(l-i+1)$ for $i\leq l${\rm;}

\item[\rm(6)] $|\omega|=|\mu|+2(p-l)+2t+1$ and $\ell(\omega)=\ell(\mu)+1$.
\end{itemize}			

Let $p_1$ and $p_2$ be the first and the second starting cluster indexes of $\mu$ respectively.
Then, we consider the following seven cases.

Case 1: $j=6$. In this case, we have $I_{p,t}(\mu)=\mu^{(2)}_{p_1}$. By the definition of $\mathcal{I}^{(6)}_{p,t}$, we see that
 $2t+1$ is marked with $1$ in $GG(\omega)$, $\omega^{(2)}_{p}=\mu^{(2)}_{p}=2t+4$, $\omega^{(2)}_{p}$ is of starting type $s_3$, and  $\omega^{(2)}_{i}=\mu^{(2)}_{i}=\mu^{(2)}_{p}+4(p-i)=\omega^{(2)}_{p}+4(p-i)$  for $p_1\leq i\leq p$. Combining with the conditions (4) and (5), we arrive at $\omega\in\mathbb{C}^{(6)}_{=}(k,r|p,t)$ and $D_{p,t}(\omega)=\omega^{(2)}_{p_1}=\mu^{(2)}_{p_1}=I_{p,t}(\mu)$.

Case 2: $j=7$. In this case, we know that $2t+1$ is marked with $1$ in $GG(\omega)$ and $\omega^{(2)}_{p+1}=2t+2$. For $p_1\leq i\leq p$, we define $r_i$ as in Definition \ref{defi-insertion-7}.
Then,  we see that  $\omega^{(2)}_i=\mu^{(2)}_i+2=\omega^{(2)}_{p+1}+4(p-i+1)$ and there is an $r_i$-marked  $\omega^{(2)}_i$ in $GG(\omega)$, where $p_1\leq i\leq p$. In particular, we have  $\omega^{(2)}_p=2t+6$. It is clear from the choice of $r_i$ that  $r_i\neq 2$,  and so $\omega^{(2)}_i$ occurs at least twice in $\omega$. In light of the conditions (4) and (5), we find that in order to prove that $\omega\in\mathbb{C}^{(7)}_{=}(k,r|p,t)$ and $D_{p,t}(\omega)=I_{p,t}(\mu)$, it remains to show that

(A)  $\omega^{(2)}_{p+1}=2t+2$ is of starting type $s_0$;

(B) $2t+2$ occurs once in $\omega$.

{\noindent Condition (A).} Suppose to the contrary that $\omega^{(2)}_{p+1}=2t+2$ is not of starting type $s_0$. Since there is a $1$-marked $2t+1$ in $GG(\omega)$, we see that $\omega^{(2)}_{p+1}=2t+2$ is of starting type $s_2$. Under the condition that $\omega^{(2)}_{p_1}=\omega^{(2)}_{p+1}+4(p-p_1+1)$, we find that  $\omega^{(2)}_{p_1}$ is of starting type $s_2$, which implies that there is a $1$-marked $\omega^{(2)}_{p_1}+2$ in  $GG(\omega)$. Note that $I_{p,t}(\mu)=\mu^{(2)}_{p_1}+2=\omega^{(2)}_{p_1}$, so we deduce that  $I_{p,t}(\mu)+2$ occurs in $\omega$, which leads to a contradiction.  Hence, $\omega^{(2)}_{p+1}=2t+2$ is of starting type $s_0$.

{\noindent Condition (B).} By the construction of $\omega$, it suffices to show that $2t+2$ occurs once in $\mu$. It follows from $\mu\in\mathbb{C}^{(7)}_{\sim}(k,r|p,t)$ that $\mu^{(2)}_p=2t+4$ is of starting type $s_2$. Then, there is no $2$-marked $2t+2$ in  $GG(\mu)$ since there is a $2$-marked $2t+4$ in $GG(\mu)$. There is also no $2$-marked $2t$ in  $GG(\mu)$, otherwise we obtain that $\mu^{(2)}_{p+1}=2t$ and it is of starting type $s_2$, which contradicts the condition (4) in Definition \ref{defi-c<}. By the definition of G\"ollnitz-Gordon marking, we find that there are no parts $2t+2$ with marks greater than $2$ in $GG(\mu)$. It  yields that $2t+2$ occurs once in $\mu$. Hence, the condition (B) holds.

Case 3: $j=8$. In this case, we find that  $\omega^{(2)}_{p+1}=2t+1$, $\omega^{(2)}_{p}=\mu^{(2)}_p+2=2t+4$ is of starting type $s_3$ and $\omega^{(2)}_{p_1}=\mu^{(2)}_{p_1}+2=\omega^{(2)}_{p}+4(p-p_1)$. It yields $I_{p,t}(\mu)=\mu^{(2)}_{p_1}+2=\omega^{(2)}_{p_1}$. Using the conditions (4) and (5), we get $\omega\in\mathbb{C}^{(8)}_{=}(k,r|p,t)$ and $D_{p,t}(\omega)=I_{p,t}(\mu)$.

Case 4: $j=9$. In this case, we see that  $2t+1$ is marked with $1$ in $GG(\omega)$,  $\omega^{(2)}_{p}=\mu^{(2)}_p=2t+2$ and $\omega^{(2)}_{p_1}=\mu^{(2)}_{p_1}=\omega^{(2)}_{p}+4(p-p_1)$,
 and so  $I_{p,t}(\mu)=\mu^{(2)}_{p_1}+2=\omega^{(2)}_{p_1}+2$.  Moreover, there are $1$-marked
$\omega^{(2)}_{p}+2=\mu^{(2)}_{p}+2,\ldots,\omega^{(2)}_{p_1}+2=\mu^{(2)}_{p_1}+2$ in $GG(\omega)$.
 It follows from $I_{p,t}(\mu)+2=\omega^{(2)}_{p_1}+4$ does not occur in $\omega$  that $\omega^{(2)}_{p},\ldots,\omega^{(2)}_{p_1}$ are of starting type $s_2$. Hence, we have
$\omega\in\mathbb{C}^{(9)}_{=}(k,r|p,t)$ and $D_{p,t}(\omega)=\omega^{(2)}_{p_1}+2=I_{p,t}(\mu)$.

Case 5: $j=10$. In this case, we obtain that  $2t+1$ is marked with $1$ in $GG(\omega)$, $\omega^{(2)}_{p+1}=\mu^{(2)}_p=2t+2$, $\omega^{(2)}_{p_1}=\mu^{(2)}_{p_1}+4=\omega^{(2)}_{p+1}+4(p-p_1+1)$, and $\omega^{(2)}_{i}=\mu^{(2)}_{i-1}=\omega^{(2)}_{p+1}+4(p-i+1)$ for $p_1< i\leq p$. In particular, we have $\omega^{(2)}_{p}=2t+6$ and $I_{p,t}(\mu)=\mu^{(2)}_{p_1}+4=\omega^{(2)}_{p_1}$.
Moreover, there are $1$-marked parts $\omega^{(2)}_{p}-2=\mu^{(2)}_{p}+2,\ldots,\omega^{(2)}_{p_1}-2=\mu^{(2)}_{p_1}+2$ in $GG(\omega)$.
Note that $I_{p,t}(\mu)+2=\omega^{(2)}_{p_1}+2$ does not occur in $\omega$, we find that $\omega^{(2)}_{p},\ldots,\omega^{(2)}_{p_1}$ are of starting type $s_1$ and $\omega^{(2)}_{p+1}=2t+2$ is of starting type $s_0$. It follows from Proposition \ref{prop-i-2}  that $\omega^{(2)}_{p_1}=\mu^{(2)}_{p_1}+4$ occurs once in $\mu$. By the construction of $\omega$, we deduce that $\omega^{(2)}_{p_1}$ occurs once in $\omega$. Appealing to the condition  (5), we get $\omega\in\mathbb{C}^{(10)}_{=}(k,r|p,t)$ and $D_{p,t}(\omega)=\omega^{(2)}_{p_1}=I_{p,t}(\mu)$.

Case 6: $j=11$.  In this case, we find that $2t+1$ is marked with $1$ in $GG(\omega)$, $\omega^{(2)}_{p}=\mu^{(2)}_{p}=2t+2$, $\omega^{(2)}_{i}=\mu^{(2)}_{i}=\omega^{(2)}_{p}+4(p-i)$ for $p_1\leq i< p$, and $\omega^{(2)}_{i}=\mu^{(2)}_{i}=\omega^{(2)}_{p}+4(p-i)+2$ for $p_2\leq i\leq p_1-1$. Moreover, there are $1$-marked parts $\omega^{(2)}_{p}+2=\mu^{(2)}_{p}+2,\ldots,\omega^{(2)}_{p_1}+2=\mu^{(2)}_{p_1}+2,\omega^{(2)}_{p_1-1}=\mu^{(2)}_{p_1},\ldots,\omega^{(2)}_{p_2}=\mu^{(2)}_{p_2}$ in $GG(\omega)$. It yields that $\omega^{(2)}_{p},\ldots,\omega^{(2)}_{p_1}$ are of starting type $s_2$ and  $\omega^{(2)}_{p_1-1},\ldots,\omega^{(2)}_{p_2}$ are of starting type $s_3$. In virtue of the conditions (4) and (5), we arrive at   $\omega\in\mathbb{C}^{(11)}_{=}(k,r|p,t)$ and $D_{p,t}(\omega)=\omega^{(2)}_{p_2}=I_{p,t}(\mu)$.

 Case 7: $j=12$.   With a similar argument as in the proofs of Case 2 and Case 5, we get  $\omega\in\mathbb{C}^{(12)}_{=}(k,r|p,t)$ and $D_{p,t}(\omega)=\omega^{(2)}_{p_2}=\mu^{(2)}_{p_2}+2$. Thus, the proof is complete.   \qed

 Now, we are in a position to define the insertion $\mathcal{I}_{p,t}$.
 \begin{defi}\label{defi-insertion}
 Let $\mu$ be a partition in $\mathbb{C}_{\sim}(k,r|p,t)$. Define $\mathcal{I}_{p,t}(\mu)=\mathcal{I}^{(j)}_{p,t}(\mu)$ if $\mu\in\mathbb{C}^{(j)}_{\sim}(k,r|p,t)$, where $1\leq j\leq 12$.
 \end{defi}

 By Lemmas \ref{lem-insertion-1-5} and \ref{lem-insertion-6-12}, we get the following lemma, which says that the insertion $\mathcal{I}_{p,t}$ is a map from $\mathbb{C}_{\sim}(k,r|p,t)$ to $\mathbb{C}_{=}(k,r|p,t)$.

 \begin{lem}\label{lem-insertion}
Let $\mu$ be a partition in $\mathbb{C}_{\sim}(k,r|p,t)$ and let $\omega=\mathcal{I}_{p,t}(\mu)$. Then, $\omega$ is a partition in $\mathbb{C}_{=}(k,r|p,t)$ such that $D_{p,t}(\omega)=I_{p,t}(\mu),$
\[|\omega|=|\mu|+2(p-l)+2t+1\text{ and }\ell(\omega)=\ell(\mu)+1,\]
where $l$ is the largest integer such that $\mu^{(2)}_{l}>I_{p,t}(\mu)$.
\end{lem}

\subsection{The separation $\mathcal{S}_{p,t}$}

In this subsection, we define the $(k-1)$-separation $\mathcal{S}_{p,t}$ from $\mathbb{C}_{=}(k,r|p,t)$ to $\mathbb{C}_{\sim}(k,r|p,t)$, which plays the role of the inverse map of the $(k-1)$-insertion $\mathcal{I}_{p,t}$.  To do this,  we will define the $j$-th kind of the separtion   $\mathcal{S}^{(j)}_{p,t}$ from $\mathbb{C}^{(j)}_{=}(k,r|p,t)$ to $\mathbb{C}^{(j)}_{\sim}(k,r|p,t)$ for $1\leq j\leq 12$.
We first present the  $j$-th kind of the separation   $\mathcal{S}^{(j)}_{p,t}$ for $1\leq j\leq 5$.
\begin{defi}\label{defi-separation-1-5}
 	
 For $1\leq j\leq 5$,	let $\omega$ be a partition in $\mathbb{C}^{(j)}_{=}(k,r|p,t)$.  Define the $j$-th kind of the separation   $\mathcal{S}^{(j)}_{p,t}$ as follows{\rm:} remove $2t+1$ from $\omega$.

 \end{defi}

We proceed to show that $\mathcal{S}^{(j)}_{p,t}$ is a map from $\mathbb{C}^{(j)}_{=}(k,r|p,t)$ to $\mathbb{C}^{(j)}_{\sim}(k,r|p,t)$ for $1\leq j\leq 5$.
\begin{lem}\label{lem-separation-1-5}
For $1\leq j\leq 5$, let $\omega$ be a partition in $\mathbb{C}^{(j)}_{=}(k,r|p,t)$ and let $\mu=\mathcal{S}^{(j)}_{p,t}(\omega)$. Then, $\mu$ is a partition in $\mathbb{C}^{(j)}_{\sim}(k,r|p,t)$ such that
\[|\mu|=|\omega|-2t-1\text{ and }\ell(\mu)=\ell(\omega)-1.\]
\end{lem}

\pf We first show that $\mu$ is a partition in $\mathbb{C}_{<}(k,r|p,t)$, that is, $\mu$ satisfies the conditions (1)-(4) in Definition \ref{defi-c<}. Note that $\mu$ is obtain by removing $2t+1$ from $\omega$, so there is no odd part of $\mu$ greater than or equal to $2t+1$, which means that $\mu$ satisfies the condition (1) in Definition \ref{defi-c<}.

By the conditions (1) and (2) in Proposition \ref{prop-divisio-index-1}, we deduce that $2t+2$ and $2t+4$ can not both occur in $\omega$, and so $2t+2$ and $2t+4$ can not both occur in $\mu$. Then, the marks of parts greater than or equal to $2t+4$ in $GG(\mu)$ are the same as those in $GG(\omega)$, which yields \begin{equation}\label{condition-(2)-1}
\mu^{(2)}_p=\omega^{(2)}_p\geq 2t+4.
\end{equation}
It implies that $\mu$ satisfies the condition (3) in Definition \ref{defi-c<}.

By the construction of $\mu$, we deduce that the marks of parts not exceeding $2t$ in $GG(\mu)$ are the same as those in $GG(\omega)$. We proceed to show that
\begin{equation}\label{label-2t+2}
2t+2\text{ can only be marked with }1\text{ in }GG(\mu).
\end{equation}
Assume that $2t+2$ occurs in $\mu$, then $2t+2$ also occurs in $\omega$. By the definition of $\mathbb{C}^{(j)}_{=}(k,r|p,t)$ for $1\leq j\leq 5$, we find that $\omega^{(2)}_{p+1}=2t+2$,  and so $\omega\in\mathbb{C}^{(1)}_{=}(k,r|p,t)$ or $\mathbb{C}^{(3)}_{=}(k,r|p,t)$. It follows from the condition (4) in Definition \ref{defi-=} that $\omega^{(2)}_{p+1}=2t+2$ is of starting type $s_0$.
Appealing to the definitions of $\mathbb{C}^{(1)}_{=}(k,r|p,t)$ and $\mathbb{C}^{(3)}_{=}(k,r|p,t)$, we see that for $i\leq p$, if
$\omega^{(2)}_i=\omega^{(2)}_{p+1}+4(p-i+1)$ then we have  $\omega\in\mathbb{C}^{(3)}_{=}(k,r|p,t)$ and $\omega^{(2)}_i$ is of starting type $s_3$.
Again by the condition (4) in Definition \ref{defi-=}, we get that $\omega^{(2)}_{p+1}=2t+2$ occurs once in $\omega$.
 Note that $\mu$ is obtained by removing the $1$-marked $2t+1$ in $GG(\omega)$, then by the definition of G\"ollnitz-Gordon marking, we see that $2t+2$ is marked with $1$ in $GG(\mu)$. Hence, \eqref{label-2t+2} is satisfied.

Using \eqref{label-2t+2}, we get $\mu^{(2)}_{p+1}<2t+2$. From the proof above, we know that $2t+1$ does not occur in $\mu$. Then, we have $\mu^{(2)}_{p+1}<2t+1$.
Combining with \eqref{condition-(2)-1}, we deduce that $\mu$ satisfies the condition (2) in Definition \ref{defi-c<}.

Assume that $\mu^{(2)}_{p+1}=2t$, then there is a $2$-marked $2t$ in $GG(\omega)$. It yields that there is no $1$-marked $2t$ in $GG(\omega)$, otherwise the mark of $2t+1$ in $GG(\omega)$ is greater than $2$, which contradicts the condition (2) in Definition \ref{defi-=}.
It follows from the definition of G\"ollnitz-Gordon marking that there is no $2$-marked $2t+2$ in $GG(\omega)$. By the condition (6) in Definition \ref{defi-=}, we find that $2t+2$ does not occur in $\omega$.
By the construction of $\mu$, we see that $2t+2$ does not occur in $\mu$ and there is no $1$-marked $2t$ in $GG(\mu)$, which implies that $\mu^{(2)}_{p+1}=2t$ is of starting type $s_0$ or $s_1$. Thus, $\mu$ satisfies the condition (4) in Definition \ref{defi-c<}.

We have proved that $\mu$ satisfies the conditions (1)-(4) in Definition \ref{defi-c<}, and so we have $\mu\in\mathbb{C}_{<}(k,r|p,t)$. Then, we wish to show  that $\mu\in\mathbb{C}^{(j)}_{\sim}(k,r|p,t)$ for $1\leq j\leq 5$.

From the proof above, we find that $2t+2$ occurs in $\mu$ if and only if $\omega^{(2)}_{p+1}=2t+2$.
By the construction of $\mu$, we obtain that  $2t+2$ does not occur in $\mu$ for $j=2$ or $4$, and
$2t+2$ occurs in $\mu$ for $j=3$.

Recall that the marks of parts greater than or equal to $2t+4$ in $GG(\mu)$ are the same as those in $GG(\omega)$, then we get

\begin{itemize}

\item[(1)] $\mu^{(2)}_p=\omega^{(2)}_p\geq 2t+8$ for $j=1$;

\item[(2)] $\mu^{(2)}_{p}=2t+6$ is of reduction type $\hat{A}_1$ or $\hat{A}_2$ or $\hat{B}$ for $j=2$;

\item[(3)] $\mu^{(2)}_{p}=2t+6$ is of reduction type $\hat{A}_1$ for $j=3$;

\item[(4)] $\mu^{(2)}_{p}=2t+6$ is of reduction type $\hat{C}$ for $j=4$;

\item[(5)] $\mu^{(2)}_{p}=2t+4$ is of reduction type $\hat{A}_1$ for $j=5$.

\end{itemize}

So far we   have accomplished the task of showing that $\mu$ is a partition in $\mathbb{C}^{(j)}_{\sim}(k,r|p,t)$ for $1\leq j\leq 5$. Recall that $\mu$ is obtain by removing $2t+1$ from $\omega$, then we have  $|\mu|=|\omega|-2t-1$  and $\ell(\mu)=\ell(\omega)-1$.  This completes the proof. \qed

Next, we give the $j$-th kind of the separation $\mathcal{S}^{(j)}_{p,t}$ for $6\leq j\leq 12$.

\begin{defi}\label{defi-separation-6}
		Let $\omega$ be a partition in $\mathbb{C}^{(6)}_{=}(k,r|p,t)$. Assume that $D_{p,t}(\omega)=\omega^{(2)}_{l+1}$. Define $\mathcal{S}^{(6)}_{p,t}\colon \omega \rightarrow \mu$ as follows{\rm:} remove  $2t+1$ from $\omega$ and replace the $1$-marked parts $\omega^{(2)}_p,\ldots,\omega^{(2)}_{l+1}$ in $GG(\omega)$ by $1$-marked parts $\omega^{(2)}_p-2,\ldots,\omega^{(2)}_{l+1}-2$ respectively to get $\mu$.
\end{defi}

For example, let $\omega$ be a partition in $\mathbb{C}^{(6)}_{=}(4,3|2,1)$ with G\"ollnitz-Gordon marking given in \eqref{example-=-6}. We have $D_{2,1}(\omega)=\omega^{(2)}_1=10$. Removing  $3$ from $\omega$ and replacing the $1$-marked parts $6$ and $10$ in $GG(\omega)$ by $1$-marked parts $4$ and $8$  respectively, we can recover the partition $\mu$ in \eqref{example-<-6}.

\begin{defi}\label{defi-separation-7}
		Let $\omega$ be a partition in $\mathbb{C}^{(7)}_{=}(k,r|p,t)$.  Assume that $D_{p,t}(\omega)=\omega^{(2)}_{l+1}$. Define $\mathcal{S}^{(7)}_{p,t}\colon \omega \rightarrow \mu$ as follows{\rm:}

	{\rm(1)} For $l+1\leq i\leq p$, let $r_i$ be the smallest integer such that $r_i\neq 2$ and there is an $r_i$-marked $\omega^{(2)}_i$ in $GG(\omega)$. Replace the $r_i$-marked  $\omega^{(2)}_i$ in $GG(\omega)$ by  $r_i$-marked $\omega^{(2)}_i-2$ for $l+1\leq i\leq p$ and denote the result partition  by $\nu$. Then, there are $1$-marked parts $\omega^{(2)}_p-2,\ldots,\omega^{(2)}_{l+1}-2$ in $GG(\nu)$.

  {\rm(2)} Remove $2t+1$ from $\nu$ to get $\mu$. Then, the parts $\omega^{(2)}_p-2,\ldots,\omega^{(2)}_{l+1}-2$ marked with $1$ in $GG(\nu)$ are marked with $2$ in $GG(\mu)$, and the parts the parts $\omega^{(2)}_{p+1},\omega^{(2)}_{p}\ldots,\omega^{(2)}_{l+1}$ marked with $2$ in  $GG(\nu)$ are marked with $1$ in $GG(\mu)$.

\end{defi}

For example, let $\omega$ be the partition in $\mathbb{C}^{(7)}_{=}(4,3|3,1)$ defined in \eqref{example-=-7}. It can be checked that  $D_{3,1}(\omega)=\omega^{(2)}_1=16$. Moreover, we have $r_3=3$, $r_2=1$ and $r_1=1$. Replacing the $3$-marked $8$, $1$-marked $12$ and $1$-marked $16$ in ${GG}(\omega)$ by $3$-marked $6$, $1$-marked $10$ and $1$-marked $14$ respectively, we obtain the partition $\nu$ in \eqref{example-<=-7}. Then, remove $3$ from $\nu$ to get the partition $\mu$ in \eqref{example-<-7}.

\begin{defi}\label{defi-separation-8}
		Let $\omega$ be a partition in $\mathbb{C}^{(8)}_{=}(k,r|p,t)$.  Assume that  $D_{p,t}(\omega)=\omega^{(2)}_{l+1}$. Define $\mathcal{S}^{(8)}_{p,t}\colon \omega \rightarrow \mu$ as follows{\rm:} remove  $2t+1$ from $\omega$ and replace the $2$-marked parts $\omega^{(2)}_p,\ldots,\omega^{(2)}_{l+1}$ in $GG(\omega)$ by $2$-marked parts $\omega^{(2)}_p-2,\ldots,\omega^{(2)}_{l+1}-2$ respectively to get $\mu$.
\end{defi}

For example, let $\omega$ be a partition in $\mathbb{C}^{(8)}_{=}(4,3|3,2)$, whose G\"ollnitz-Gordon marking is given in \eqref{example-=-8}. We have $D_{3,2}(\omega)=\omega^{(2)}_1=16$. Removing  $5$ from $\omega$ and replacing the $2$-marked parts $8$, $12$ and $16$ in $GG(\omega)$ by $2$-marked parts $6$, $10$ and $14$ respectively, we can recover the partition $\mu$ in \eqref{example-<-8}.

\begin{defi}\label{defi-separation-9}
Let $\omega$ be a partition in $\mathbb{C}^{(9)}_{=}(k,r|p,t)$.  Assume that  $D_{p,t}(\omega)=\omega^{(2)}_{l+1}+2$. Define $\mathcal{S}^{(9)}_{p,t}\colon \omega \rightarrow \mu$ as follows{\rm:} remove  $2t+1$ from $\omega$ and replace the $1$-marked parts $\omega^{(2)}_p+2,\ldots,\omega^{(2)}_{l+1}+2$ in $GG(\omega)$ by $1$-marked parts $\omega^{(2)}_p,\ldots,\omega^{(2)}_{l+1}$ respectively to get $\mu$.
\end{defi}

For example, let $\omega$ be a partition in $\mathbb{C}^{(9)}_{=}(4,3|3,0)$ given in \eqref{example-=-9}. We have $D_{3,0}(\omega)=\omega^{(2)}_1+2=12$. Remove  $1$ from $\omega$ and replace the $1$-marked parts $4$, $8$ and $12$ in $GG(\omega)$ by $1$-marked parts $2$, $6$ and $10$ respectively to recover the partition $\mu$ in \eqref{example-<-9}.

\begin{defi}\label{defi-separation-10}
		Let $\omega$ be a partition in $\mathbb{C}^{(10)}_{=}(k,r|p,t)$.  Assume that  $D_{p,t}(\omega)=\omega^{(2)}_{l+1}$. Define $\mathcal{S}^{(10)}_{p,t}\colon \omega \rightarrow \mu$ as follows{\rm:} remove $2t+1$ from $\omega$ and replace the $1$-marked parts $\omega^{(2)}_{p+1}+2,\ldots,\omega^{(2)}_{l+2}+2$ in $GG(\omega)$ by $1$-marked parts $\omega^{(2)}_{p+1},\ldots,\omega^{(2)}_{l+2}$ respectively to get $\mu$. Then, the  part $\omega^{(2)}_{l+1}$ marked with $2$ in  $GG(\omega)$ is  marked with   $1$ in $GG(\mu)$.
\end{defi}

For example, let $\omega$ be a partition in $\mathbb{C}^{(10)}_{\sim}(4,3|3,0)$ with G\"ollnitz-Gordon marking given in \eqref{example-=-10}. Clearly, $D_{3,0}(\omega)=\omega^{(2)}_1=14$.
Remove  $1$ from $\omega$ and replace the $1$-marked parts $4$, $8$ and $12$ in $GG(\omega)$ by $1$-marked parts $2$, $6$ and $10$ respectively to get the partition $\mu$ in \eqref{example-<-10}.
Then, the  part $14$ marked with $2$ in  $GG(\omega)$ is  marked with   $1$ in $GG(\mu)$.

\begin{defi}\label{defi-separation-11}
Let $\omega$ be a partition in $\mathbb{C}^{(11)}_{=}(k,r|p,t)$. Assume that $D_{p,t}(\omega)=\omega^{(2)}_{l+1}$ and $s$ is the smallest integer such that $\omega^{(2)}_{s}=\omega^{(2)}_{p}+4(p-s)$.
Define $\mathcal{S}^{(11)}_{p,t}\colon \omega \rightarrow \mu$ as follows{\rm:}  remove $2t+1$ from $\omega$ and  replace the $1$-marked parts $\omega^{(2)}_p+2,\ldots,\omega^{(2)}_{s}+2,\omega^{(2)}_{s-1},\ldots,\omega^{(2)}_{l+1}$ in $GG(\omega)$ by $1$-marked parts $\omega^{(2)}_p,\ldots,\omega^{(2)}_{s},\omega^{(2)}_{s-1}-2,\ldots,\omega^{(2)}_{l+1}-2$ respectively to get $\mu$.
\end{defi}

For example, let $\omega$ be a partition in $\mathbb{C}^{(11)}_{=}(4,3|5,0)$, whose G\"ollnitz-Gordon marking is given in \eqref{example-=-11}. It can be checked that $D_{5,0}(\omega)=\omega^{(2)}_1=20$ and $s=3$. Removing  $1$ from $\omega$ and replacing the $1$-marked $4$, $8$, $12$, $16$ and $20$ in $GG(\omega)$ by $1$-marked $2$, $6$, $10$, $14$ and $18$  respectively, we recover the partition in \eqref{example-<-11}.

\begin{defi}\label{defi-separation-12}
								
Let $\omega$ be a partition in $\mathbb{C}^{(12)}_{=}(k,r|p,t)$. Assume that $D_{p,t}(\omega)=\omega^{(2)}_{l+1}$ and $s$ is the smallest integer such that $\omega^{(2)}_{s}=\omega^{(2)}_{p}+4(p-s)$ and $\omega^{(2)}_{s}$ occurs once in $\omega$. Define $\mathcal{S}^{(12)}_{p,t}\colon \omega \rightarrow \mu$ as follows{\rm:}

{\rm(1)} For $l+1\leq i< s$, let $r_i$ be the smallest integer such that $r_i\neq 2$ and there is an $r_i$-marked  $\omega^{(2)}_i$ in $GG(\omega)$.  Replace the $r_i$-marked  $\omega^{(2)}_i$ in $GG(\omega)$ by  $r_i$-marked part $\omega^{(2)}_i-2$ for $l+1\leq i< s$ and denote the result partition  by $\nu$.
Then, there are $1$-marked parts $\omega^{(2)}_{s-1}-2,\ldots,\omega^{(2)}_{l+1}-2$ in $GG(\nu)$.

{\rm(2)}  Remove $2t+1$ from $\nu$ and  replace the $1$-marked parts $\omega^{(2)}_{p+1}+2,\ldots,\omega^{(2)}_{s+1}+2$ in $GG(\nu)$ by $1$-marked parts $\omega^{(2)}_{p+1},\ldots,\omega^{(2)}_{s+1}$ respectively to get $\mu$.
Then, the parts $\omega^{(2)}_{s-1}-2,\ldots,\omega^{(2)}_{l+1}-2$ marked with $1$ in  $GG(\nu)$ are marked with $2$ in $GG(\mu)$, and the parts $\omega^{(2)}_{s},\omega^{(2)}_{s-1}\ldots,\omega^{(2)}_{l+1}$ marked with $2$ in  $GG(\nu)$ are marked with $1$ in $GG(\mu)$.
								
\end{defi}

For example, let $\omega$ be a partition in $\mathbb{C}^{(12)}_{=}(4,3|6,0)$ defined in \eqref{example-=-12}. It can be checked that $D_{6,0}(\omega)=\omega^{(2)}_1=26$, $s=4$, $r_3=3$, $r_2=1$ and $r_1=1$.  Replacing the $3$-marked $18$, $1$-marked $22$ and $1$-marked $26$ in $GG(\omega)$ by $3$-marked $16$, $1$-marked $20$ and $1$-marked $24$ respectively, we get the partition $\nu$ in \eqref{example-<=-12}. Then, remove $1$ from $\nu$ and replace the $1$-marked $4$, $8$ and $12$ in $GG(\nu)$ by $1$-marked $2$, $6$ and $10$ respectively to recover the partition in \eqref{example-<-12}.

We proceed to show that for $6\leq j\leq 12$, the $j$-th kind of the separation $\mathcal{S}^{(j)}_{p,t}$ is a map from $\mathbb{C}^{(j)}_{=}(k,r|p,t)$ to $\mathbb{C}^{(j)}_{\sim}(k,r|p,t)$.	

\begin{lem}\label{lem-separation-6-12}
For $6\leq j\leq 12$, let $\omega$ be a partition in $\mathbb{C}^{(j)}_{=}(k,r|p,t)$ and let $\mu=\mathcal{S}^{(j)}_{p,t}(\omega)$. Then, $\mu$ is a partition in $\mathbb{C}^{(j)}_{\sim}(k,r|p,t)$ such that $I_{p,t}(\mu)=D_{p,t}(\omega),$
\[|\mu|=|\omega|-2(p-l)-2t-1\text{ and }\ell(\mu)=\ell(\omega)-1,\]
where $l$ is the largest integer such that $\omega^{(2)}_{l}>D_{p,t}(\omega)$.
\end{lem}

\pf It can be checked that the marks of the unchanged parts not exceeding  $D_{p,t}(\omega)$ in $GG(\mu)$ are the same as those in $GG(\omega)$. It is clear from the conditions (3) in Proposition \ref{prop-divisio-index-1} that $D_{p,t}(\omega)+2$ does not occur in $\omega$. By the construction of $\mu$, we see that $D_{p,t}(\omega)+2$ does not occur in $\mu$. So, the marks of parts greater than  $D_{p,t}(\omega)+2$ in $GG(\mu)$ are the same as those in $GG(\omega)$. Hence, the $j$-th kind of the separation $\mathcal{S}^{(j)}_{p,t}$ is well-defined for $6\leq j\leq 12$.

By the choice of $D_{p,t}(\omega)$, we find that for $i\leq l$, if  $\omega^{(2)}_i=D_{p,t}(\omega)+4(l-i+1)$ then $\omega^{(2)}_i$ is of starting type $s_3$ and $\omega^{(2)}_i+2$ occurs in $\omega$. By the construction of $\mu$, we obtain that for $i\leq l$, if $\mu^{(2)}_i=D_{p,t}(\omega)+4(l-i+1)$ then $\mu^{(2)}_i$ is of starting type $s_3$ and $\mu^{(2)}_i+2$ occurs in $\mu$. For $6\leq j\leq 12$, in order to prove that $\mu\in\mathbb{C}^{(j)}_{\sim}(k,r|p,t)$, it suffices to show that $\mu\in\mathbb{C}^{(j)}_{<}(k,r|p,t)$ and $I_{p,t}(\mu)=D_{p,t}(\omega)$.
Clearly, there is no odd part of $\mu$ greater than or equal to $2t+1$. Then, we consider the following seven cases.

Case 1: $j=6$. In this case, we have $\mu^{(2)}_p=\omega^{(2)}_p=2t+4$, $\mu^{(2)}_{p+1}=\omega^{(2)}_{p+1}<2t+1$ and $D_{p,t}(\omega)=\omega^{(2)}_{l+1}=\mu^{(2)}_{l+1}$. Moreover, there are $1$-marked parts $\mu^{(2)}_p-2,\ldots,\mu^{(2)}_{l+1}-2$ in $GG(\mu)$ and $\mu^{(2)}_{l+1}+2$ does not occur in $\mu$, and so $\mu^{(2)}_p,\ldots,\mu^{(2)}_{l+1}$ are of starting type $s_1$. It implies that $\mu\in\mathbb{C}^{(6)}_{<}(k,r|p,t)$ and $I_{p,t}(\mu)=\mu^{(2)}_{l+1}=D_{p,t}(\omega)$.

Case 2: $j=7$. In this case, we have $\mu^{(2)}_p=\omega^{(2)}_p-2=2t+4$, $\mu^{(2)}_{p+1}=\omega^{(2)}_{p+2}<2t$ and $D_{p,t}(\omega)=\omega^{(2)}_{l+1}=\mu^{(2)}_{l+1}+2$. Moreover, there are $1$-marked parts $\mu^{(2)}_p+2,\ldots,\mu^{(2)}_{l+1}+2$ in $GG(\mu)$ and $\mu^{(2)}_{l+1}+4=\omega^{(2)}_{l+1}+2$ does not occur in $\mu$, and so $\mu^{(2)}_p,\ldots,\mu^{(2)}_{l+1}$ are of starting type $s_2$. It implies that $\mu\in\mathbb{C}^{(7)}_{<}(k,r|p,t)$ and $I_{p,t}(\mu)=\mu^{(2)}_{l+1}+2=D_{p,t}(\omega)$.

Case 3: $j=8$. In this case, we have $\mu^{(2)}_p=\omega^{(2)}_p-2=2t+2$, $\mu^{(2)}_{p+1}=\omega^{(2)}_{p+2}<2t$ and $D_{p,t}(\omega)=\omega^{(2)}_{l+1}=\mu^{(2)}_{l+1}+2$. With a similar argument as in Case 2, we arrive at $\mu\in\mathbb{C}^{(8)}_{<}(k,r|p,t)$ and $I_{p,t}(\mu)=\mu^{(2)}_{l+1}+2=D_{p,t}(\omega)$.

Case 4: $j=9$. In this case, we have $\mu^{(2)}_p=\omega^{(2)}_p=2t+2$, $\mu^{(2)}_{p+1}=\omega^{(2)}_{p+1}<2t$ and $D_{p,t}(\omega)=\omega^{(2)}_{l+1}+2=\mu^{(2)}_{l+1}+2$. Moreover, there are $1$-marked parts $\mu^{(2)}_p,\ldots,\mu^{(2)}_{l+1}$ in $GG(\mu)$ and $\mu^{(2)}_{l+1}+4=\omega^{(2)}_{l+1}+4$ does not occur in $\mu$. It yields that  $\mu^{(2)}_p,\ldots,\mu^{(2)}_{l+1}$ are of starting type $s_3$, $\mu\in\mathbb{C}^{(9)}_{<}(k,r|p,t)$ and $I_{p,t}(\mu)=\mu^{(2)}_{l+1}+2=D_{p,t}(\omega)$.

Case 5: $j=10$. In this case, we have $\mu^{(2)}_p=\omega^{(2)}_{p+1}=2t+2$, $\mu^{(2)}_{p+1}=\omega^{(2)}_{p+2}<2t$ and $D_{p,t}(\omega)=\omega^{(2)}_{l+1}=\mu^{(2)}_{l+1}+4$. Moreover, there are $1$-marked parts $\mu^{(2)}_p,\ldots,\mu^{(2)}_{l+1}$ in $GG(\mu)$, and so  $\mu^{(2)}_p,\ldots,\mu^{(2)}_{l+1}$ are of starting type $s_3$. Since $\omega\in\mathbb{C}^{(10)}_{=}(k,r|p,t)$, we see that $D_{p,t}(\omega)=\omega^{(2)}_{l+1}$ occurs once in $\omega$ and $D_{p,t}(\omega)+2=\omega^{(2)}_{l+1}+2$ does not occur in $\omega$. It follows from the construction of $\mu$ that $D_{p,t}(\omega)=\omega^{(2)}_{l+1}=\mu^{(2)}_{l+1}+4$ occurs once in $\mu$ and $D_{p,t}(\omega)+2=\omega^{(2)}_{l+1}+2=\mu^{(2)}_{l+1}+6$ does not occur in $\omega$. Hence, we have $\mu\in\mathbb{C}^{(10)}_{<}(k,r|p,t)$ and $I_{p,t}(\mu)=\mu^{(2)}_{l+1}+4=D_{p,t}(\omega)$.

Case 6: $j=11$. In this case, we have $\mu^{(2)}_p=\omega^{(2)}_{p}=2t+2$, $\mu^{(2)}_{p+1}=\omega^{(2)}_{p+1}<2t$ and $D_{p,t}(\omega)=\omega^{(2)}_{l+1}=\mu^{(2)}_{l+1}$.
Assume $s$ is the smallest integer such that $\omega^{(2)}_{s}=\omega^{(2)}_{p}+4(p-s)$, then there are $1$-marked parts $\mu^{(2)}_p,\ldots,\mu^{(2)}_{s},\mu^{(2)}_{s-1}-2,\ldots,\mu^{(2)}_{l+1}-2$ in $GG(\mu)$. It implies that
 $\mu^{(2)}_p,\ldots,\mu^{(2)}_{s}$ are of starting type $s_3$.  It follows from the construction of $\mu$ that $D_{p,t}(\omega)+2=\mu^{(2)}_{l+1}+2$ does not occur in $\mu$, and so $\mu^{(2)}_{s-1},\ldots,\mu^{(2)}_{l+1}$ are of starting type $s_1$.
Note that $\mu^{(2)}_{s-1}=\omega^{(2)}_{s-1}=\omega^{(2)}_{s}+6=\mu^{(2)}_{s}+6$, we obtain that $\mu\in\mathbb{C}^{(11)}_{<}(k,r|p,t)$ and $I_{p,t}(\mu)=\mu^{(2)}_{l+1}=D_{p,t}(\omega)$.

Case 6: $j=12$. With a similar argument as in Case 2 and Case 5,  we get $\mu\in\mathbb{C}^{(12)}_{<}(k,r|p,t)$ and $I_{p,t}(\mu)=\mu^{(2)}_{l+1}+2=\omega^{(2)}_{l+1}=D_{p,t}(\omega)$.

Now, we  conclude that for $6\leq j\leq 12$, $\mu$ is a partition in $\mathbb{C}^{(j)}_{<}(k,r|p,t)$ and $I_{p,t}(\mu)=D_{p,t}(\omega)$, and so $\mu\in\mathbb{C}^{(j)}_{\sim}(k,r|p,t)$. Evidently,  $|\mu|=|\omega|-2(p-l)-2t-1$  and $\ell(\mu)=\ell(\omega)-1$.  This completes the proof. \qed

 Now, we are in a position to define the separation $\mathcal{S}_{p,t}$.
 \begin{defi}\label{defi-separation}
 Let $\omega$ be a partition in $\mathbb{C}_{=}(k,r|p,t)$. Define $\mathcal{S}_{p,t}(\omega)=\mathcal{S}^{(j)}_{p,t}(\omega)$ if $\omega\in\mathbb{C}^{(j)}_{=}(k,r|p,t)$, where $1\leq j\leq 12$.
 \end{defi}

 By Lemmas \ref{lem-separation-1-5} and \ref{lem-separation-6-12}, we get the following lemma, which says that the separation $\mathcal{S}_{p,t}$ is a map from $\mathbb{C}_{=}(k,r|p,t)$ to $\mathbb{C}_{\sim}(k,r|p,t)$.

 \begin{lem}\label{lem-separation}
Let $\omega$ be a partition in $\mathbb{C}_{=}(k,r|p,t)$ and let $\mu=\mathcal{S}_{p,t}(\omega)$. Then, $\mu$ is a partition in $\mathbb{C}_{\sim}(k,r|p,t)$ such that $I_{p,t}(\mu)=D_{p,t}(\omega),$
\[|\mu|=|\omega|-2(p-l)-2t-1\text{ and }\ell(\mu)=\ell(\omega)-1,\]
where $l$ is the largest integer such that $\omega^{(2)}_{l}>D_{p,t}(\omega)$.
\end{lem}

\subsection{Proof of Theorem \ref{equiv-main-I}}

We conclude this section with a proof of Theorem \ref{equiv-main-I}.

{\bf \noindent Proof of Theorem \ref{equiv-main-I}:} Using lemma \ref{lem-insertion}, we obtain that the insertion $\mathcal{I}_{p,t}$ is a map from $\mathbb{C}_{\sim}(k,r|p,t)$ to $\mathbb{C}_{=}(k,r|p,t)$ satisfying \eqref{equiv-main-I-0-eqn}.
Appealing to lemma \ref{lem-separation}, we get that the separation $\mathcal{S}_{p,t}$ is a map from $\mathbb{C}_{=}(k,r|p,t)$ to $\mathbb{C}_{\sim}(k,r|p,t)$. It follows from the definitions of  $\mathcal{I}_{p,t}$ and  $\mathcal{S}_{p,t}$ that they are inverse of each other. This completes the proof. \qed

\section{Proof of Theorem \ref{main-thm}}

 We first give an equivalent combinatorial statement of Theorem \ref{main-thm}. Let $\mathbb{E}(k,r)$ denote the number of partitions in $\mathbb{C}(k,r)$ without odd parts.
 For $N_1\geq N_2\geq \cdots\geq N_{k-1}\geq 0$, let $\mathbb{E}(N_1,\ldots,N_{k-1};r)$ denote the set of partitions $\pi$ in $\mathbb{E}(k,r)$ with $N_i$ parts in the $i$-th row of $GG(\pi)$ for $1\leq i\leq k-1$.
 Based on Gordon marking, Kur\c{s}ung\"{o}z \cite{Kursungoz-2010a,Kursungoz-2010} established the following identity.
\begin{equation*}
\sum_{\pi\in\mathbb{E}(N_1,\ldots,N_{k-1};r)}x^{\ell(\pi)}q^{\frac{|\pi|}{2}}=\frac{x^{N_1+\cdots+N_{k-1}}q^{N_1^2+\cdots+N_{k-1}^2+N_r+\cdots+N_{k-1}}}{(q;q)_{N_1-N_2}\cdots(q;q)_{N_{k-2}-N_{k-1}} (q;q)_{N_{k-1}}}.
\end{equation*}
Then, we can see that
\begin{equation*}
\sum_{\pi\in\mathbb{E}(N_1,\ldots,N_{k-1};r)}q^{|\pi|}=\frac{x^{N_1+\cdots+N_{k-1}}q^{2(N_1^2+\cdots+N_{k-1}^2+N_r+\cdots+N_{k-1})}}{(q^2;q^2)_{N_1-N_2}\cdots(q^2;q^2)_{N_{k-2}-N_{k-1}} (q^{2};q^{2})_{N_{k-1}}}.
\end{equation*}

For $N\geq 0$, let $\mathbb{I}_{N}$ denote the set of partitions $\zeta=(2m_1+1,2m_2+1,\ldots,2m_\ell+1)$ with distinct odd parts greater than or equal to $2N+1$, that is, $m_1>m_2>\cdots>m_{\ell}\geq N$.  The generating function for partitions in $\mathbb{I}_{N}$ is
\[\sum_{\zeta\in\mathbb{I}_{N}}x^{\ell(\zeta)}q^{|\zeta|}=(1+xq^{2N+1})(1+xq^{2N+3})\cdots=(-xq^{2N+1};q^2)_\infty.\]

We define $\mathbb{F}(3,3)$ to be set of pairs $(\pi,\zeta)$  of partitions such that
\[\pi\in\mathbb{E}(3,3)\text{ and }\zeta\in\mathbb{I}_{N_2(\pi)}.\]
Then, Theorem \ref{main-thm} is equivalent to the following combinatorial statement.
\begin{thm}\label{equiv-main-thm-0000-2}
There is a bijection $\Phi$ between $\mathbb{F}(3,3)$ and $\mathbb{C}(3,3)$. Moreover, for a pair $(\pi,\zeta)\in \mathbb{F}(3,3)$, we have $\omega=\Phi(\pi,\zeta)\in\mathbb{C}(3,3)$ such that
\[|\omega|=|\pi|+|\zeta| \text{ and } \ell(\omega)=\ell(\pi)+\ell(\zeta).\]
\end{thm}

To give a proof of Theorem \ref{equiv-main-thm-0000-2}, we need to introduce two sets $\mathbb{C}_{<}(k,r|m)$ and $\mathbb{C}_{=}(k,r|m)$, and build a bijection $\Phi_{m}$ between $\mathbb{C}_{<}(k,r|m)$ and $\mathbb{C}_{=}(k,r|m)$ based on the bijection $\Phi_{p,t}=\mathcal{I}_{p,t}\cdot\mathcal{H}_{p,t}$. Then, we can construct the bijection $\Phi$ in Theorem \ref{equiv-main-thm-0000-2} with the aid of the bijection $\Phi_{m}$.

\subsection{$\mathbb{C}_{<}(k,r|m)$ and $\mathbb{C}_{=}(k,r|m)$}

 For $m\geq 0$, we set
\[\mathbb{C}_{<}(k,r|m)=\bigcup_{p+t=m}\mathbb{C}_{<}(k,r|p,t).\]
The following proposition states that for $m\geq 0$ and $\pi\in\mathbb{C}_{<}(k,r|m)$, there exist unique integers $p$ and $t$ such that $p+t=m$ and $\pi\in\mathbb{C}_{<}(k,r|p,t)$.

\begin{prop}\label{unique-<}
For $m\geq 0$, let $\pi$ be a partition in $\mathbb{C}_{<}(k,r|p,t)$, where $p+t=m$. Then, there do not exist  integers $p'$ and $t'$ such that $p'\neq p$, $t'\neq  t$,  $p'+t'=m$ and $\pi\in \mathbb{C}_{<}(k,r|p',t')$.
\end{prop}

\pf Suppose to the contrary that there exist integers $p'$ and $t'$ such that $p'\neq p$, $t'\neq  t$,  $p'+t'=m$ and $\pi\in \mathbb{C}_{<}(k,r|p',t')$. Without loss of generality, we assume that $p'<p$. By definition, we have
\begin{equation}\label{lem-pro-unique-1}
2t+1<\pi^{(2)}_p\leq\pi^{(2)}_{p'+1}<2t'+1,
\end{equation}
 which yields
\begin{equation}\label{prop-unique-2}
2(t'-t)>\pi^{(2)}_{p'+1}-\pi^{(2)}_p.
\end{equation}

It follows from the definition of G\"ollnitz-Gordon marking that
\begin{equation*}\label{prop-unique-3}
\pi^{(2)}_{p'+1}-\pi^{(2)}_p\geq 4(p-p'-1).
\end{equation*}
Combining with \eqref{prop-unique-2}, we get
\[2(t'-t)>4(p-p'-1).\]

Note that $p+t=p'+t'=m$, so we obtain that $p'>p-2$. Under the assumption that $p'<p$, we have $p'=p-1$, and so $t'=t+1$. Using \eqref{lem-pro-unique-1}, we find that
 \[\pi^{(2)}_p=2t+2\text{ and }\pi^{(2)}_{p'+1}=2t'.\]

  Since $\pi$ is a partition in $\mathbb{C}_{<}(k,r|p,t)$, we see that  $\pi^{(2)}_p=2t+2$ is of starting type $s_2$ or $s_3$. But, under the assumption that $\pi\in\mathbb{C}_{<}(k,r|p',t')$, we obtain that $\pi^{(2)}_{p'+1}=2t'$ is of starting type $s_0$ or $s_1$, which leads to a contradiction. Thus, we complete the proof.  \qed

The following theorem gives a   criterion   to determine whether a partition in $\mathbb{E}(k,r)$ is also a partition in $\mathbb{C}_{<}(k,r|m)$.

\begin{thm}\label{einc<}

For $N_2\geq 0$, let $\pi$ be a partition in $\mathbb{E}(k,r)$ such that there are $N_2$ parts marked with $2$ in $GG(\pi)$. Then, $\pi$ is a partition in $\mathbb{C}_{<}(k,r|m)$ if and only if $m\geq N_2$.

\end{thm}

\pf We first show that if $m\geq N_2$, then $\pi$ is a partition in $\mathbb{C}_{<}(k,r|m)$. Assume that $l$ is the largest integer such that $2(m-l)+1<\pi^{(2)}_l$. Such an integer $l$ exists because $\pi^{(2)}_0=+\infty$. By the choice of $l$, we get  $\pi^{(2)}_{l+1}<2(m-l)$.

It is clear from $\pi\in\mathbb{E}(k,r)$ that there are no odd parts in $\pi$. If $\pi^{(2)}_l=2(m-l)+2$ and it is of starting type $s_1$, then we have $\pi\in \mathbb{C}_{<}(k,r|l-1,m-l+1)$. Otherwise, we have $\pi\in \mathbb{C}_{<}(k,r|l,m-l)$. In either case, we can get $\pi\in \mathbb{C}_{<}(k,r|m)$. This completes the proof of the sufficiency.

 Conversely, assume that $\pi$ is a partition in $\mathbb{C}_{<}(k,r|m)$, then by Proposition \ref{geqN2}, we get $m\geq N_2$. This completes the proof.   \qed

 For $m\geq 0$, we set
\[\mathbb{C}_{=}(k,r|m)=\bigcup_{p+t=m}\mathbb{C}_{=}(k,r|p,t).\]
By definition, we have
\begin{prop}\label{unique-=}
For $m\geq 0$, let $\pi$ be a partition in $\mathbb{C}_{=}(k,r|m)$. Then, there exist unique integers $p$ and $t$ such that $p+t=m$ and $\pi\in\mathbb{C}_{=}(k,r|p,t)$.
\end{prop}

We consider the case $k=r=3$.

\begin{thm}\label{33e=}
Let $\pi$ be a partition in $\mathbb{C}(3,3)$ such that there exist odd parts in $\pi$. Then,
there exists unique $m$ such that $\pi\in\mathbb{C}_{=}(3,3|m)$.
\end{thm}

\pf Assume that the largest odd part in $\pi$ is $2t+1$. Let $l$ be the largest integer such that $\pi^{(2)}_{l}> 2t+1$. In light of Lemma \ref{prop-new-add} and the condition (1) in Corollary \ref{2t+2}, we obtain that if $2t+2$ occurs in $\pi$ then $2t+2$ can only be marked with $2$ in $GG(\pi)$ and it is of starting type $s_0$ or $s_2$.
Then, we consider the following two cases.

Case 1: $\pi^{(2)}_{l}>2t+2$, or $\pi^{(2)}_{l}=2t+2$ with starting type $s_2$. Obviously, we have $\pi\in\mathbb{C}_{=}(3,3|l,t)$, and so $\pi\in\mathbb{C}_{=}(3,3|l+t)$.

Case 2: $\pi^{(2)}_{l}=2t+2$ and it is of starting type $s_0$. Under the condition that $2t+2$ can only be marked with $2$ in $GG(\pi)$, we know that $2t+2$ occurs once in $\pi$. It yields $\pi\in\mathbb{C}_{=}(3,3|l-1,t)$, and so $\pi\in\mathbb{C}_{=}(3,3|l+t-1)$.

In either case, we have shown that $\pi$ is a partition in  $\mathbb{C}_{=}(3,3|m)$ for an unique $m$. The proof is complete.   \qed

For $m\geq 0$, let $\pi$ be a partition in  $\mathbb{C}_{<}(k,r|m)$, define
\[\Phi_m(\pi)=\Phi_{p,t}(\pi)\text{ if }\pi\in\mathbb{C}_{<}(k,r|p,t).\]
By Theorem \ref{equiv-main-lem-1-1}, Proposition \ref{unique-<} and Proposition \ref{unique-=}, we get the following theorem.
\begin{thm}\label{equiv-main-thm-2}
For $m\geq0,$ the map $\Phi_{m}$ is a bijection between $\mathbb{C}_{<}(k,r|m)$ and $\mathbb{C}_{=}(k,r|m)$. Moreover, for a partition $\pi\in \mathbb{C}_{<}(k,r|m)$, we have $\omega=\Phi_{m}(\pi)\in\mathbb{C}_{=}(k,r|m)$ such that
\[|\omega|=|\pi|+2m+1 \text{ and } \ell(\omega)=\ell(\pi)+1.\]
\end{thm}

It is worth mentioning that the inverse map $\Psi_{m}$ of $\Phi_{m}$ is defined as follows. For $m\geq 0$, let $\omega$ be a partition in  $\mathbb{C}_{=}(k,r|m)$, define
\[\Psi_m(\omega)=\mathcal{R}_{p,t}(\mathcal{S}_{p,t}(\omega))\text{ if }\omega\in\mathbb{C}_{=}(k,r|p,t).\]
The following lemma is an immediate consequence of  Lemma \ref{equiv-main-R-lem-1}, Lemma \ref{lem-separation}, Proposition \ref{unique-<} and Proposition \ref{unique-=}.
\begin{lem}\label{equiv-main-lemm-000}
For $m\geq0,$ the map $\Psi_{m}$ is a map from $\mathbb{C}_{=}(k,r|m)$ to $\mathbb{C}_{<}(k,r|m)$. Moreover, for a partition $\omega\in \mathbb{C}_{=}(k,r|m)$, we have $\pi=\Psi_{m}(\omega)\in\mathbb{C}_{<}(k,r|m)$ such that
\[|\pi|=|\omega|-2m-1 \text{ and } \ell(\pi)=\ell(\omega)-1.\]
\end{lem}

We conclude this subsection with the following theorem, which involves  the successive application of $\Phi_m$.

\begin{thm}\label{ssins}

For $m\geq 0$, let $\pi$ be a partition in $\mathbb{C}_{=}(k,r|m)$. Then, $\pi$ is a partition in $\mathbb{C}_{<}(k,r|m')$ if and only if $m<m'$.

\end{thm}

\pf Since $\pi$ be a partition in $\mathbb{C}_{=}(k,r|m)$, there exist unique $p$ and $t$ such that  $p+t=m$ and $\pi\in\mathbb{C}_{=}(k,r|p,t)$.  By definition, we have $\pi^{(2)}_{p}\geq 2t+2$, $\pi^{(2)}_{p+1}\leq 2t+2$, and the largest odd part of $\pi$ is $2t+1$.

We first show that if $m<m'$ then  $\pi$ is  in $\mathbb{C}_{<}(k,r|m')$. Assume that $l$ is the largest integer such that $2(m'-l)+1<\pi^{(2)}_l$. Recall that $\pi^{(2)}_{p+1}\leq 2t+2$, so we have   $\pi^{(2)}_{p+2}\leq 2t-1$. Under the condition that $m<m'$, we get
\[2(m'-p-2)+1\geq 2(m-p-1)+1=2t-1\geq \pi^{(2)}_{p+2},\]
 which yields $l\leq p+1$, and so
 \begin{equation}\label{use-0000}
 \pi^{(2)}_l>2(m'-l)+1\geq2(m+1-p-1)+1=2t+1.
  \end{equation}
  Note that the largest odd part of $\pi$ is $2t+1$, we see that $\pi^{(2)}_l$ is not of type $s_{-1}$.


 If $\pi^{(2)}_l=2(m'-l)+2$ and it is of type $s_0$ or $s_1$, then we set $p'=l-1$. Otherwise, we set $p'=l$.  Let $t'=m'-p'$.  We proceed to show that $\pi\in\mathbb{C}_{<}(k,r|p',t')$, which yields $\pi\in\mathbb{C}_{<}(k,r|m')$. To do this, it remains to prove that the largest odd part of $\pi$ is less than $2t'+1$, namely, $t<t'$. Under the assumption that $m<m'$, we just need to show that $p'\leq p$. Suppose to the contrary that $p'=l=p+1$. Combining with \eqref{use-0000}, we have $\pi^{(2)}_{p+1}=\pi^{(2)}_l>2t+1$.
 Recall that $\pi^{(2)}_{p+1}\leq 2t+2$,  we find that $\pi^{(2)}_l=\pi^{(2)}_{p+1}=2t+2$. Again by \eqref{use-0000}, we get  $m'=m+1$. So, we have $\pi^{(2)}_{p+1}=2t+2=2(m'-l)+2$.
 From the definition of $\mathbb{C}_{=}(k,r|p,t)$,  we deduce that $\pi^{(2)}_{p+1}=2t+2$ is of type $s_0$. From the proof above, we have  $p'=l-1$, which contradicts the assumption that $p'=l$. Hence, we have shown $p'\leq p$. This completes the proof of the sufficiency.

 Conversely, assume that $\pi$ is a partition in $\mathbb{C}_{<}(k,r|m')$, we intend to show that $m<m'$. Suppose to the contrary that $m\geq m'$. Assume that $\pi\in\mathbb{C}_{<}(k,r|p',t')$, where   $p'+t'=m'$.  Recall that the largest odd part of $\pi$ is $2t+1$, so we have $t<t'$. Then we have $p>p'$. Hence, we get
 \[\begin{split}
 2m'+1&=2(p'+t')+1=2t'+1+2p'\\
 &>\pi^{(2)}_{p'+1}+2p'\geq\pi^{(2)}_{p}+2(p-p'-1)+2p'\\
 &\geq 2t+2+2(p-1)=2m,
 \end{split}\]
which implies that $m'\geq m$. Under the assumption that  $m\geq m'$, we have $m=m'$. Moreover, we obtain that $p'=p-1$, $t'=t+1$ and $\pi^{(2)}_{p}=2t+2$. Under the assumption that $\pi\in\mathbb{C}_{<}(k,r|p',t')$, we see that   $\pi^{(2)}_{p'+1}=\pi^{(2)}_{p}=2t'$ is of type $s_0$ or $s_1$. But, $\pi^{(2)}_{p'+1}=\pi^{(2)}_{p}=2t+2$ is of type $s_2$ since $\pi\in\mathbb{C}_{=}(k,r|p,t)$, which leads to a contradiction. Thus, we have shown $m<m'$. This completes the proof.   \qed

\subsection{Proof of Theorem \ref{equiv-main-thm-0000-2}}

We are now in a position to give a proof of Theorem \ref{equiv-main-thm-0000-2}.

 {\noindent \bf Proof of Theorem \ref{equiv-main-thm-0000-2}.} Let $(\pi,\zeta)$ be a pair in $\mathbb{F}(3,3)$. Then, we have $\pi\in\mathbb{E}(3,3)$. We define $\omega=\Phi(\pi,\zeta)$ as follows. There are two cases.

 Case 1: $\zeta=\emptyset$. Then we set $\omega=\pi$. It is clear that $\omega\in\mathbb{E}(3,3)\subseteq \mathbb{C}(3,3)$. Moreover, $|\omega|=|\pi|+|\zeta|$ and $\ell(\omega)=\ell(\pi)+\ell(\zeta)$.

 Case 2: $\zeta\neq\emptyset$.  Assume that there are $N_2$ parts marked with $2$ in $GG(\pi)$. Then, we have $\zeta\in\mathbb{I}_{N_2}$, denoted  $\zeta=(2m_1+1,2m_2+1,\ldots,2m_\ell+1)$, where $m_1>m_2>\cdots>m_\ell\geq N_2$.
 Starting with $\pi$, we apply the bijection $\Phi_m$ repeatedly to get $\omega$. Denote the intermediate partitions by $\pi^0,\pi^1,\ldots,\pi^\ell$ with $\pi^0=\pi$ and $\pi^\ell=\omega$. By Theorem \ref{einc<}, we have  $\pi^0=\pi\in\mathbb{C}_<(3, 3|m_\ell)$.

Set $b=0$ and repeat the following process until $b=\ell$.

\begin{itemize}
    \item[(A)] Note that $\pi^b\in\mathbb{C}_<(3,3|m_{\ell-b})$, we apply $\Phi_{m_{\ell-b}}$ to $\pi^b$ to get $\pi^{b+1}$, that is,
   \[\pi^{b+1}=\Phi_{m_{\ell-b}}(\pi^{b}).\]
   By Theorem \ref{equiv-main-thm-2}, we deduce that
   \[
   \pi^{b+1}\in\mathbb{C}_=(3,3|m_{\ell-b}),\]
   \begin{equation*}\label{phisum}
   |\pi^{b+1}|=|\pi^b|+2m_{\ell-b}+1,
   \end{equation*}
   and
   \begin{equation*}\label{philength}
   \ell(\pi^{b+1})=\ell(\pi^b)+1.
   \end{equation*}

    \item[(B)] Replace $b$ by $b+1$. If $b=\ell$, then we are done. If $b<\ell$, then by Theorem \ref{ssins}, we have
     \[
   \pi^{b}\in\mathbb{C}_<(3,3|m_{\ell-b}),\]
 since $m_{\ell-b}>m_{\ell-b+1}$. Go back to (A).
\end{itemize}

Eventually, the above process yields $\omega=\pi^\ell\in \mathbb{C}_=(3,3|m_{1})$ such that $|\omega|=|\pi|+|\zeta|$ and $\ell(\omega)=\ell(\pi)+\ell(\zeta)$. Moreover, we have $\omega\in\mathbb{C}(3,3)$.

 To show that $\Phi$ is a bijection, we give the inverse map $\Psi$ of $\Phi$ by successively applying $\Psi_m$. Let $\pi$ be a partition in $\mathbb{C}(3,3)$. We shall construct a pair $(\pi,\zeta)$, that is, $(\pi,\zeta)=\Psi(\omega)$, such that $(\pi,\zeta)\in \mathbb{F}(3,3)$, $|\omega|=|\pi|+|\zeta|$ and $\ell(\omega)=\ell(\pi)+\ell(\zeta)$. Assume that there are $\ell\geq 0$ odd parts in $\omega$. We eliminate all odd parts of $\omega$. There are two cases.

 Case 1: $\ell=0$. Then set $\pi=\omega$ and $\zeta=\emptyset$. Clearly, $(\pi,\zeta)\in \mathbb{F}(3,3)$, $|\omega|=|\pi|+|\zeta|$ and $\ell(\omega)=\ell(\pi)+\ell(\zeta)$.

 Case 2: $\ell\geq 1$.  We eliminate the $\ell$ odd parts of $\pi$  by successively applying $\Psi_m$. Denote the intermediate pairs by $(\omega^0,\zeta^0),(\omega^1,\zeta^1)\ldots,(\omega^\ell,\zeta^\ell)$ with  $(\omega^0,\zeta^0)=(\omega,\emptyset)$.

 Set $b=0$ and carry out the following procedure.

 \begin{itemize}
     \item[(A)] Since there exist odd parts in $\omega^{b}$, then by Theorem \ref{33e=}, we see that there exists $m_{b+1}$ such that
      \[
   \omega^{b}\in\mathbb{C}_=(3,3|m_{b+1}).\]
     Apply $\Psi_{m_{b+1}}$ to $\omega^b$ to get $\omega^{b+1}$, that is,
      \[\omega^{b+1}=\Psi_{m_{b+1}}(\omega^{b}).\]
      Employing Lemma \ref{equiv-main-lemm-000}, we find that
   \[
   \omega^{b+1}\in\mathbb{C}_<(3,3|m_{b+1}),\]
    \begin{equation*}\label{psisum}
   |\omega^{b+1}|=|\omega^b|-2m_{b+1}-1,
   \end{equation*}
   and
   \begin{equation*}\label{psilength}
   \ell(\omega^{b+1})=\ell(\omega^b)-1.
   \end{equation*}
     Then insert $2m_{b+1}+1$ into $\zeta^b$ to obtain $\zeta^{b+1}$.

     \item[(B)] Replace $b$ by $b+1$. If $b=\ell$, then we are done. Otherwise, go back to (A).

 \end{itemize}

Observe that for $0\leq b\leq \ell$, there are $\ell-b$ odd parts in $\omega^b$. In particular, there are no odd parts in $\omega^\ell$, and so $\omega^\ell\in\mathbb{E}(3,3)$. Assume that there are $N_2$ parts marked with $2$ in $GG(\omega^\ell)$. Note that $\omega^{\ell}\in\mathbb{C}_<(3,3|m_{\ell})$, then by Theorem \ref{33e=}, we get
\[m_\ell\geq N_2.\]
Theorem \ref{ssins} reveals that for $0\leq b<\ell-1$,
\[m_{b+1}>m_{b+2}.\]
Therefore, we see that $\zeta^\ell=(2m_1+1,2m_2+1,\ldots,2m_{\ell}+1)$ is a partition in $\mathbb{I}_{N_2}$.
Set $(\pi,\zeta)=(\omega^\ell,\zeta^\ell)$.  It is clear that $\pi=\omega^\ell\in\mathbb{E}(3,3)$, $\zeta=\zeta^\ell\in \mathbb{I}_{N_2}$, and so  $(\pi,\zeta)$ is a pair in  $\mathbb{F}(3,3)$. Moreover, we have $|\omega|=|\pi|+|\zeta|$ and $\ell(\omega)=\ell(\pi)+\ell(\zeta)$.

Note that $\Psi_m$ is the inverse of $\Phi_m$, it can be verified that $\Psi$ is the inverse of $\Phi$. Thus, we complete the proof. \qed


\begin{thebibliography}{99}

 \bibitem{Andrews-1976} G.E. Andrews, The Theory of partitions, Addison-Wesley Publishing Co., 1976.

\bibitem{Bressoud-1980} D.M. Bressoud, Analytic and combinatorial generalizations of Rogers-Ramanujan identities, Mem. Amer. Math. Soc. 24 (227) (1980) 54pp.

    \bibitem{Gollnitz-1960}H. G\"ollnitz, Einfache Partionen, Diplomarbeit W. S. Gottingen, 65 , 1960.

         \bibitem{Gollnitz-1967}H. G\"ollnitz, Partitionen mit differenzenbedingungen, J. reine angew. Math. 225, 1967, 154--190.

             \bibitem{Gordon-1961} B. Gordon, A combinatorial generalization of the Rogers-Ramanujan identities, Amer. J. Math 83 (1961) 393--399.


\bibitem{Gordon-1962} B. Gordon, Some Ramanujan-like continued fractions. Abstracts of Short Communications. Int. 1962.



    \bibitem{Gordon-1965} B. Gordon, Some continued fractions of the Rogers-Ramanujan type, Duke Math. J. 31, 1965, 741--748.

   \bibitem{He-Zhao-2023} Thomas Y. He and Alice X.H. Zhao, New companions to the generations of the G\"ollnitz-Gordon identities, Ramanujan J. 61 (2023) 1077--1120.


 \bibitem{Kim-2018} S. Kim, Bressoud's conjecture, Adv. Math.     325 (2018) 770--813.

   \bibitem{Kim-Yee-2014} S. Kim and A.J. Yee, Partitions with part difference conditions and Bressoud's conjecture, J. Comb. Theory Ser. A 126 (2014) 35--69.

   \bibitem{Kursungoz-2010a} K. Kur\c{s}ung\"{o}z, Parity considerations in Andrews-Gordon identities and the $k$-marked Durfee symbols, PhD thesis, Penn State University, 2009.

\bibitem{Kursungoz-2010} K. Kur\c{s}ung\"{o}z, Parity considerations in Andrews-Gordon identities, Eur. J. Comb. 31 (2010) 976--1000.

\end{thebibliography}
\end{document}